\documentclass[11pt,fullpage]{article}
\usepackage{chngcntr}

\usepackage{amsmath}
\usepackage{dsfont}
\usepackage{mathrsfs}
\usepackage{graphicx}
\usepackage{amsfonts}
\usepackage{amssymb}
\usepackage[round]{natbib}
\usepackage{bbding} 
\usepackage{bm}
\usepackage{enumitem}
\usepackage{hhline}
\usepackage{comment}
\usepackage{multirow}
\usepackage{caption}
\usepackage{subcaption}
\usepackage{float}
\restylefloat{table}
\usepackage[title]{appendix}%

\usepackage[english,activeacute]{babel}
\usepackage{graphicx}
\usepackage{amssymb}

\usepackage[usenames]{color}
\usepackage[pagebackref=true, colorlinks=true, citecolor=blue, anchorcolor=blue, urlcolor=blue]{hyperref}
\usepackage[hang,flushmargin]{footmisc}
\usepackage{amsthm}
\usepackage{my_definitions}
\usepackage[capitalise,noabbrev]{cleveref}

\usepackage[mathlines]{lineno}[4.41]
  \newcommand*\patchAmsMathEnvironmentForLineno[1]{%
    \expandafter\let\csname old#1\expandafter\endcsname\csname #1\endcsname
    \expandafter\let\csname oldend#1\expandafter\endcsname\csname end#1\endcsname
    \renewenvironment{#1}%
    {\linenomath\csname old#1\endcsname}%
    {\csname oldend#1\endcsname\endlinenomath}}%
  \newcommand*\patchBothAmsMathEnvironmentsForLineno[1]{%
    \patchAmsMathEnvironmentForLineno{#1}%
    \patchAmsMathEnvironmentForLineno{#1*}}%
  \patchBothAmsMathEnvironmentsForLineno{equation}%
  \patchBothAmsMathEnvironmentsForLineno{align}%
  \patchBothAmsMathEnvironmentsForLineno{flalign}%
  \patchBothAmsMathEnvironmentsForLineno{alignat}%
  \patchBothAmsMathEnvironmentsForLineno{gather}%
  \patchBothAmsMathEnvironmentsForLineno{multline}

\parskip\medskipamount

\newcommand{\eps}{\epsilon}              
\newcommand{\ind}{1 \!\!1}              

\newcommand{\comps}{{\sf comps}}

\newtheorem{prop}{Proposition}

\newcommand{\bi}{\begin{itemize}}
\newcommand{\ei}{\end{itemize}}
\newcommand{\be}{\begin{enumerate}}
\newcommand{\ee}{\end{enumerate}}

\newtheorem{assum}{Assumption}
\newtheorem{thm}{Theorem}

\newtheorem{dfn}{Definition}
\newtheorem{cor}{Corollary}
\newtheorem{lem}{Lemma}

\newcounter{counterproblem}
\DefineNamedColor{named}{OliveGreen}    {cmyk}{0.64,0,0.95,0.40}

\DefineNamedColor{named}{BrickRed}      {cmyk}{0,0.89,0.94,0.28}

\newcommand*\samethanks[1][\value{footnote}]{\footnotemark[#1]}

\title{Heavy Traffic Analysis of Multi-Class Bipartite Queueing Systems Under FCFS\\
}

\author{{\large Lisa Aoki Hillas}\thanks{Corresponding author. Booth School of Business, The University of Chicago. Email: lhillas@chicagobooth.edu}\hspace{1.5cm} Ren\'e Caldentey\thanks{Booth School of Business, The University of Chicago. Email: \{rene.caldentey,varun.gupta \}@chicagobooth.edu}
\hspace{1.5cm}Varun Gupta\samethanks 
}

\begin{document}
\maketitle
\parindent 0em

\begin{abstract} 
This paper examines the performance of multi-class multi-server bipartite queueing systems under a FCFS-ALIS service discipline, where each arriving customer is only compatible with a subset of servers. We analyze the system under conventional heavy-traffic conditions, where the traffic intensity approaches one from below. Building upon the formulation and results of \cite{Afecheetal2019}, we generalize the model by allowing the vector of arrival rates to approach the heavy-traffic limit from an arbitrary direction. We characterize the steady-state waiting times of the various customer classes and demonstrate that a much wider range of waiting time outcomes is achievable. Furthermore, we establish that the matching probabilities, i.e., the probabilities of different customer classes being served by different servers, do not depend on the direction along which the system approaches heavy traffic.
We also investigate the design of compatibility between customer classes and servers, finding that a service provider who has complete control over the matching can design a delay-minimizing menu by considering only the limiting arrival rates. When some constraints on the compatibility structure exist, the direction of convergence to heavy-traffic affects which menu minimizes delay. Additionally, we discover that the bipartite matching queueing system exhibits a form of Braess's paradox, where adding more connectivity to an existing system can lead to higher average waiting times, despite the fact that neither customers nor servers act strategically.


\end{abstract}
\setcounter{footnote}{1}
\section{Introduction}\label{sec:introduction}
In this paper, we analyse the performance of multi-class bipartite queuing systems under an FCFS-ALIS service discipline. Multi-class bipartite queuing systems are used for modelling a variety of important applications, such as public housing, health-care, and manufacturing. However, these models can be both analytically and computationally intractable, making questions of performance analysis and system design difficult to answer. Heavy-traffic scaling can be used to provide approximations of these systems that are much simpler to analyse and reveal fundamental properties of the system.

The specific model we look at has $n$ customer classes and $m$ distinct servers. Customers arrive to each class according to independent Poisson processes. Service times are exponentially distributed, with service rates depending only on the server, and not on the customer class. Each customer class has a particular subset of servers they can be served by. Each server may potentially be compatible with multiple customer classes. Servers serve the customer classes they are compatible with according to a FCFS discipline. That is, when a server finishes serving a customer, they consider all of the customers that belong to classes they are compatible with, and serve the customer that has been waiting the longest. 

We analyse two aspects of the performance of this model, the expected waiting time delays of the different customer classes, and the matching probabilities of the different customer classes, that is, the probability with which a customer of a given class is served by a particular server. 

This paper is an extension of \cite{Afecheetal2019}, who study a similar model to ours. Their model uses a specific heavy-traffic scaling, which limits the range of outcomes that the model can produce. By using a more general heavy-traffic scaling, we increase the range of outcomes the model produces, allowing for accurate approximations of a wider range of scenarios. Additionally, we allow for some queues to have no arrivals at the heavy-traffic limit, and are able to calculate the expected delay should a customer join such a queue. Our main motivation for considering this generalisation is the study of systems with strategic customers, i.e., customers who can choose their class type upon arrival based on waiting time delays and matching probabilities. For example, consider a system with two independent M/M/1 queues, both being served at rate 1. Further suppose that arriving customers would prefer to be served by a particular server, but also incur some waiting cost. If customers are able to choose which queue to join, and make their decisions by trading off the cost of waiting against the value of being served by the preferred server, then we would expect the average waiting time at the queue served by the preferred server to be higher than the average waiting time at the queue served by the less preferred server. 

Using a conventional heavy traffic scaling, in which the number of servers and the service rates remain fixed, and the traffic intensity approaches 1 from below, the limiting arrival rates of both queues will be 1. The heavy traffic scaling in \cite{Afecheetal2019} has the proportion of customers arriving into the different queues remaining constant while taking the limit. However, if we do this in our simple M/M/1 example, we can see that this would limit us to concluding that the heavy traffic delays of both queues are equal. However, if we generalise the approach to heavy-traffic, allowing the arrival rates into the different queues to approach their limits at different rates, we are able to increase the range of outcomes we can model. We can interpret the different rates of approach in the real world as the different queues having arrival rates closer or further away to their predicted limiting value. 

This application of strategic arrivals also motivates us to allow for queues with no arrivals. This can be important for developing a coherent model when including strategic behaviour. In this case, it is possible to offer queues that no customers will choose to join, but we still need to calculate expected delays for those queues in order to justify why customers are not choosing to join them. 

In this paper, we calculate the expected delays of the different customer classes using our more general scaling. We show different approaches to the heavy-traffic limit produce different waiting time outcomes. In \cref{sec:discussion}, we use this to show that very minor perturbations in arrival rates can produce significant improvements in waiting time outcomes in the pre-limit. Additionally, we show that the limiting matching probabilities do not depend on the scaling used, but only depend on the limiting arrival rates. Finally, we look at some simple questions regarding the design of the compatibility between customer classes and servers. We find that when the service provider has complete control over the compatibility structure, they only need to consider the limiting arrival rates in order to design a delay minimising menu. When there are some constraints on the compatibility structure, then the particular approach to heavy-traffic does affect which menu minimises delay. 

\vspace{0.3cm}
{\bf Related Literature.} Heavy-traffic approximations have long been used to simplify the study of intractable queueing systems. Early works in this area include \cite{kingman1962queues} and \cite{whitt1974heavy}. These papers look at a so-called ``conventional'' approach to heavy-traffic, in which the number of servers and their service capacities remain fixed, and the arrival rate grows large in such a way that the traffic intensity of the system converges to one from below. An alternative class of ``many-server'' heavy-traffic limits have also been considered in the literature by carefully letting the number of servers and arrival rate grow unboundedly, e.g., \cite{halfin1981heavy} or \cite{Atar2012}. Motivated by mathematical tractability as well as by the fact that many real-world service systems operate under high levels of congestion\footnote{For example, the Chicago Housing Authority reported more than 170,000 families waiting for public housing in 2021. Similarly, in the same year, about 113,589 children in the United States were waiting to be adopted. In the healthcare system, more than 100,000 people are waiting for an organ transplant at any given moment in time, with average waiting times that can be as long as 5 years for a kidney transplant according to the National Kidney Foundation.}, we will study the performance of our multi-class multi-server bipartite queuing system operating under conventional heavy traffic conditions.

A range of questions can be answered using heavy-traffic approximations. In the context of parallel service systems, \cite{Harrison1999} study the question of optimal control of parallel service systems, that is, which servers should be used to serve which customer classes, and in which order should the different customer classes be served. \cite{Harrison1999} solve an approximating Brownian control problem, and conjecture that a discrete review policy will minimise holding costs for the original queuing system. This approach of using an approximating Brownian control problem to develop an optimal policy was originally suggested by  \cite{harrison1988brownian}. \cite{williams2000dynamic} and \cite{bell2001dynamic} go on to prove the asymptotic optimality of a continuous review policy for a two-server system. Following this work, \cite{mandelbaum2004scheduling} proves the asymptotic optimality of the $c\mu-$rule for convex holding costs. A distinctive feature in all of these papers is that they impose a {\em complete resource pooling} condition on the connectivity and/or compatibility between customer classes and servers (see \citealp{Harrison1999}). Roughly speaking, this condition boils down to assuming that the servers' capacities can be pooled together so that the servers can essentially act as a single ``super-server''. This assumption significantly simplifies the analysis as it allows us to obtain a single-dimensional state-space description of the workload of the system in the heavy traffic limit.

The complete resource pooling assumption is quite restrictive, however, and can be shown not to hold when strategic customer behaviour is allowed as in \cite{caldenteyetal2023}. There has already been some work moving beyond the complete resource pooling assumption. \cite{kushner2000optimal} prove the convergence to the heavy-traffic limit of a particular class of systems that do not satisfy the complete resource pooling assumption under quite general conditions. \cite{pesic2016dynamic} generalises \cite{Harrison1999} beyond the complete resource pooling assumption.  Other works analysing multi-class multi-server queueing systems with no complete resource pooling assumption include \cite{shahAsymptoticIndependenceServers2016} and \cite{hurtado2022heavy}. \cite{shahAsymptoticIndependenceServers2016} look at a system in which servers simultaneously work to process the same job, while \cite{hurtado2022heavy} analyse a generalised switch problem under a MaxWeight service policy. 

In addition to studying the problem of optimal control, questions regarding the performance of parallel service systems have been studied using heavy-traffic approximations, or fluid approximations more generally. \cite{Talreja2008} looks at the problem of calculating matching rates for a parallel service system operating under FCFS, that is, with what probability is each customer class served by each server, although the authors looked at this question for an overloaded system with abandonments. Matching rates were calculated for specific classes of networks. Various approximation methods have been developed for calculating matching rates including 
the {\em dissipative} algorithm proposed by \cite{CaldenteyKaplan2002}, a related approximation based on Ohm's law proposed by \cite{Fazel-Kaplan2018} and a quadratic programming formulation proposed by  \cite{Afecheetal2019}. 
Of these papers looking at the performance of parallel service systems under FCFS, \cite{Afecheetal2019} is the only one to also look at calculating waiting times as we do here.  Another contribution of \cite{Afecheetal2019} is to study the question of the design of matching topologies fixing the scheduling policy. While \cite{Afecheetal2019} studies this design question for a FCFS service discipline, \cite{varma2021transportation} studies the same question of the design of matching topologies under a MaxWeight service discipline. 

The specific model we look at here is a generalisation of \cite{Afecheetal2019}, which itself developed out of a long history of papers studying bipartite queueing systems and bipartite matching models under an FCFS service discipline. Early papers in this area include \cite{Schwartz1974} and \cite{Green1985}, who look at the steady-state performance of these systems given a particular hierarchical compatibility structure between customer classes and service classes, and \cite{Kaplan1984,Kaplan1988}, who similarly analysed the steady-state performance of parallel queuing systems, but for more general compatibility structures. Following \cite{Kaplan1984,Kaplan1988}, Kaplan's multi-class multi-server queueing model was adapted by \cite{CaldenteyKaplan2002}, who introduced an infinite-bipartite matching model to analyse matching probabilities under a FCFS service discipline. The model of \cite{CaldenteyKaplan2002} was further developed by \cite{CaldenteyKaplanWeiss} and then adapted by \cite{AdanWeiss2014} to that of a multi-class multi-server parallel queuing system, which is the model we use here. 

Since the development of the infinite matching model and the queueing model, different authors have looked at different aspects of the problem. \cite{BusicGuptaMairesse2013}, \cite{Mairesse2016}, and \cite{Moyal2017} look at stability conditions of such systems, and find that the system will be stable so long as a set of Hall's type conditions are satisfied. Also of interest are the steady-state matching probabilities. \cite{CaldenteyKaplanWeiss} were able to use a particular Markov chain representation to calculate the steady-state distribution of the matching system for particular classes of matching topologies. \cite{AdanWeiss2012} came up with an alternative Markov chain representation to derive the steady-state distribution of the matching system for general matching topologies, while \cite{AdanWeiss2014} used a similar approach to look at the multi-class multi-server queueing problem, and showed the equivalence of the steady-state outcomes for the matching and the overloaded queueing system. However, the combinatorial structure of the state space description of the Markov chain limits the size of the systems that can be studied both analytically and computationally. \cite{Afecheetal2019} use heavy traffic analysis to unveil a number of structural properties embedded in the infinite matching model and its corresponding multi-class bipartite matching queueing system (see also the survey by \citealp{GardnerRighter20} for a comprehensive review of related papers and models).

The rest of the paper is organized as follows. In \cref{sec:Model} we provide a detailed mathematical description of the bipartite queueing model, review some related results in the literature and introduce the heavy traffic regime that we will use to analyze the performance of the system. \cref{sec:heavytrafficwaits} is devoted to the derivation of the limiting steady-state waiting times of the different service classes. Our main result in this section is \cref{thm:CRPdelay_append} which provides a complete characterization of these limiting waiting times in terms of an underlying set of complete resource pooling components and their connectivity that emerge under heavy traffic. In \cref{sec:matchingprob} we study the steady-state matching probabilities between customer classes and servers and show in \cref{thm:matching_independent_gamma} that these probabilities do not depend on the particular direction along which the system reaches heavy traffic. This is in direct contrast to the behaviour of the steady-state waiting times, which are particularly sensitive to the direction of convergence. In \cref{sec:discussion} we discuss a number of insights that emerge from our theoretical results. For instance, what vectors of delays are implementable and how to design the connectivity between service classes and servers to achieve them. We also show that adding more connectivity to an existing bipartite queueing system can lead to longer average delays (i.e.,  some form of Braess's paradox). \cref{sec:ProofTheorem} contains the proofs and additional discussion of our main results Theorems \ref{thm:CRPdelay_append} and \ref{thm:matching_independent_gamma}. Some concluding remarks and possible directions in which our work can be extended are present in \cref{sec:Conlusion}. Finally, the Appendix contains additional proofs of various intermediate results.

\section{Model Description}\label{sec:Model}
In this section, we provide a detailed mathematical description of the model and basic definitions. To simplify our notation, we will adopt the following conventions throughout the paper. For a positive integer $k$, $[k]:=\{1,2,\dots,k\}$. All vectors are column vectors, and for a  vector $x \in \mathbb{R}^k$, $|x|:=\sum_{i \in [k]} x_i$.\vspace{0.1cm}

We consider a service system as follows. We have a set of $m$ servers organised into a set of $n$ customer classes. Each customer class is served by a particular subset of servers. This information is encoded in a compatibility matrix $M \in \{0, 1\}^{n \times m}$, where customer class $i$ can be served by server $j$ iff $m_{ij} = 1$. Customers arrive to the customer classes according to independent Poisson processes. We let $\lambda = (\lambda_1, ..., \lambda_n)$ be the arrival rates into the different customer classes. Service times are exponentially distributed, and depend only on the server. The vector of service rates will be denoted by $\mu = (\mu_1, ..., \mu_m)$. Servers will serve customers they are compatible with according to a FCFS-ALIS service discipline.  

\begin{figure}[hbt]
    \begin{center}
    \includegraphics[width=8cm]{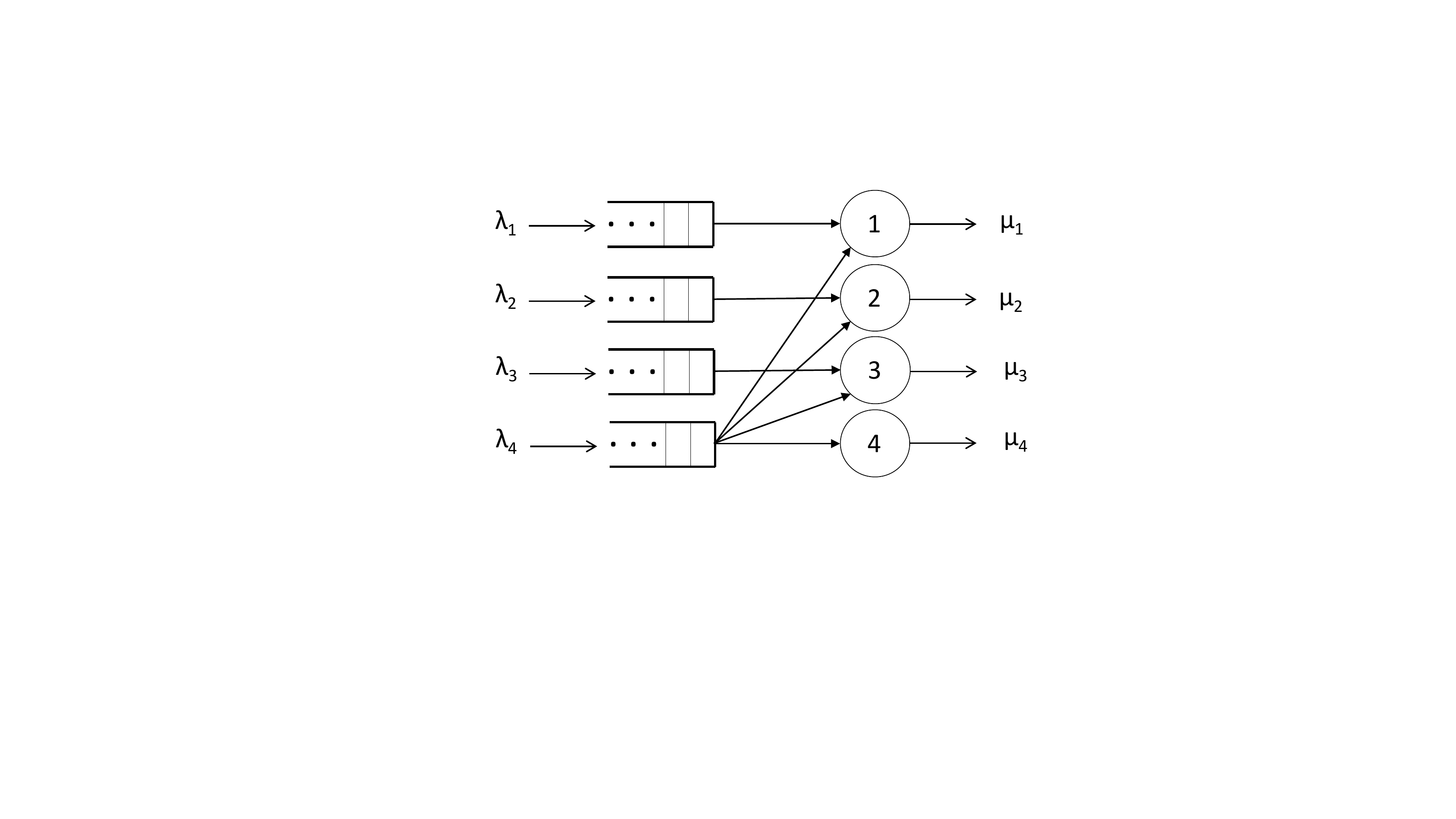}
    \caption[Model Description]{\sf \footnotesize Example with four service classes and four servers.}
     \label{fig:Ex1}
    \end{center}
\end{figure}

To illustrate, Figure~\ref{fig:Ex1} depicts an example with four servers ($m=4$), and four service classes ($n=4$). 

In this example, the menu $M$ is given by
\begin{equation}
M=\left[\begin{array}{cccc} 1 & 0 & 0 & 0 \\ 0 & 1 & 0 & 0 \\ 0 & 0 & 1 & 0 \\ 1 & 1 & 1 & 1\end{array}\right],
\end{equation}
that is, class 1 is compatible with server 1; class 2 is compatible with server 2; class 3 is compatible with server 3; and class 4 is compatible with all servers. Note that a server may belong to multiple service classes. 

We are only interested in systems which operate with stable queue lengths. The following result, from \cite{AdanWeiss2014} tells us exactly which triplets $(\lambda, \mu, M)$ produce stable steady-state outcomes. 
\begin{prop}\label{prop:stability}{\rm \citep[Theorem 2.1]{AdanWeiss2014}}
For a menu $M$ with arrival rates $\lambda$  and service rates $\mu$, define the slack of a set of servers $\Delta_{\mathscr{S}} \subseteq [m]$ as 
\begin{equation}
\Delta_{\mathscr{S}}(M) := \sum_{j \in \mathscr{S}} \mu_j - \sum_{i \in U_{\mathscr{S}}(M)} \lambda_i\qquad \mbox{for all } \mathscr{S} \subseteq [m],
\end{equation}
where
\begin{equation*}
U_{\mathscr{S}}(M):= \Big\{i \in [n]\colon \sum_{j \in \mathscr{S}^c} m_{ij}=0\Big\}
\end{equation*}
is the subset of service classes that can only be served by servers in $\mathscr{S}$. 

The menu M admit a steady state under a FCFS-ALIS service discipline if and only if:
\begin{equation*}
    \Delta_{\mathscr{S}}(M) > 0 \qquad \mbox{for all } \mathscr{S} \subseteq [m].
\end{equation*}
\end{prop}

\subsection{Steady state results for fixed arrival rates}\label{sec:adanweissresult}
Our results build on the steady state analysis of \cite{AdanWeiss2014}, which we briefly review for completeness. The authors derive their results based on a  Markov chain representation of the system defined on a carefully crafted state space $X$. A state in this state space is described by three components: ({\it i}) a permutation of servers $s = (s_1, \ldots, s_m)$, ({\it ii}) an integer $b \in \{0,\ldots, m\}$ indicating the number of busy servers, and ({\it iii}) a vector $(n_1, \ldots, n_b)$ that indicates the composition of customers waiting for service in the different service classes. It is helpful to denote a generic state $x \in X$ by the tuple:
\begin{equation}\label{eq:generic_state}
x = (s_1, n_1, s_2, n_2, \ldots, s_b, n_b, s_{b+1}, \ldots , s_m).    
\end{equation}

The first $b$ components $(s_1,\dots,s_b)$ of the server permutation $s$ denote the $b$ busy servers ranked according to the arrival time of the customer they are serving, with server $s_1$ serving the oldest arrival and server $s_b$ serving the youngest arrival. The remaining servers $(s_{b+1}, \ldots, s_m)$ are all idle and ranked in the order they became idle, with $s_{b+1}$ the server that has been idle the longest. Finally, $n_\ell$ for $\ell=1,\dots,b$, represents the number of customers in the system who arrived after the job currently being served by $s_{\ell}$ but before the job currently being served by $s_{\ell + 1}$. Due to the FCFS-ALIS service discipline, we know these customers can only be served by some server in $(s_1,\dots,s_{\ell})$. That is, each of these $n_\ell$ customers must belong to some service menu in $U( s_1,\ldots, s_\ell )$.

According to \cite[Theorem 2.1]{AdanWeiss2014}, the steady-state probability of state $x$ admits the product form:
\begin{align}\label{eq:steadystate}
\pi(x) &= {\cal B}  \prod_{\ell=1}^b \frac{\lambda^{n_\ell}_{U(s_1,\ldots, s_\ell)} }{\mu_{\{s_1,\ldots, s_\ell\}}^{n_\ell+1} } \prod_{\ell=b+1}^m \lambda_{C(s_\ell, \ldots, s_m)}^{-1},
\end{align}
where ${\cal B}$ is an appropriate normalizing constant. Additionally, each of the $n_\ell$ customers `between' server $s_\ell$ and server $s_{\ell+1}$ belongs to service class $i \in U(s_1,\ldots, s_\ell)$ independently with probability $\frac{\lambda_i}{\lambda_{U(s_1,\ldots, s_\ell)}}$.

These steady-state probabilities can be used to calculate the expected number of customers of each type in the system. Little's Law can then be applied to calculate expected steady-state mean waiting times. However, if we consider the process for calculating expected waiting times even for our relatively simple example in \cref{fig:Ex1}, we see that while these calculations are possible, the process is laborious and the resulting expressions are unwieldy. For example, let us consider how we would calculate the expected number of class 4 customers. We first observe that class 4 customers are compatible with all servers. This means that the only times class 4 customers are waiting in the system is if all servers are busy when a class 4 customer arrives. Thus if we want to calculate the expected number of class 4 customers waiting for service in the system, we can restrict ourselves to considering only the states in which all 4 servers are busy. 

Fixing the permutation of servers, and the number of busy servers, the values of $n_i$ are geometrically distributed, and hence the expected values have closed form expressions. For example, if we condition on being in the subset of states $x \in X_{(s_1,s_2,s_3,s_4)}$ such that $b = 4$ and the server permutation $(s_1, s_2, s_3, s_4)$, i.e. $x = (s_1, n_1, s_2, n_2, s_3, n_3, s_4, n_4)$, then the expected value of $n_4$ is 
\begin{equation}
\mathbb{E}(n_4| x \in X_{(s_1,s_2,s_3,s_4)}) =  \frac{\cal{B} |\lambda| \cdot |\mu|}{(\mu_1 - \lambda_1)(\mu_1 + \mu_2 - (\lambda_1 + \lambda_2))(|\mu|-\mu_4  - (|\lambda| - \lambda_4 ))(|\mu| - |\lambda|)}
\end{equation}
where $|\lambda|:=\lambda_1 + \lambda_2 + \lambda_3 + \lambda_4$, $|\mu|:=\mu_1 + \mu_2 + \mu_3 + \mu_4$ and $\cal{B}$ is an appropriate normalizing constant. Note that $n_4$ is not the number of class 4 customers; instead $n_4$ is the number of customers who arrived to the system after the customer server 4 is currently serving. Therefore the expected number of class 4 customers conditional on being in the subset of states $X_{(s_1,s_2,s_3,s_4)}$ is $\frac{\lambda_4}{|\lambda|}\mathbb{E}[n_4|x \in X_{(s_1,s_2,s_3,s_4)}]$.

To fully calculate the expected number of class 4 customers, we would need to repeat this process for every permutation of servers. Since there are four servers, there are 24 possible permutations of servers to sum over, with different combinations of terms appearing in the denominator for each permutation. This gives us very complicated expressions for the expected number of servers. If we were instead looking at the number of class 1 customers, we would also need to consider states in which only some servers are busy, giving us even more server combinations that we need to consider. 

It is this underlying computational complexity -which grows combinatorially fast in the size of the system-- that motivates our move to heavy traffic. As the system approaches heavy traffic, the probability of being in a state with an idle server approaches 0, letting us restrict our attention only to states in which all servers are busy. Additionally, we show in \cref{prop:steadystate} that in heavy-traffic, only certain server permutations have positive probability, which is a fact that simplifies the problem even further. 

\subsection{Heavy traffic scaling}
The last part of the model is the heavy-traffic scaling. As mentioned in the Introduction, our formulation extends \cite{Afecheetal2019}, who consider a specific direction of convergence to heavy traffic to derive their results. Specifically, they assume that the proportions of customers of different types remain constant as the system approaches heavy traffic. In this paper, we allow a general direction of convergence. 

We consider a conventional heavy traffic regime in which the arrival rates approach the capacity of the service system, while the number of customer classes and servers, and the service menu remain constant. We parameterize our systems by $\epsilon$, and let the service system approach heavy traffic as $\eps \downarrow 0$. Specifically, we assume that there is a sequence of arrival rates $\htp{\lambda} = \{ \htp{\lambda_i} \}_{i \in [n]}$ where 
\begin{align}
\label{eq:customer_arrival_rate}
    \htp{\lambda_i} = \Lambda_i - \gamma_i \epsilon + o(\epsilon) \geq 0\quad \mbox{for all } i \in [n] \text{ and }0 < \epsilon < \epsilon_+,
\end{align}
for some vector $\Lambda \in \RR_{+}^n$, some vector $\gamma \in \RR^n$, and some $\epsilon_+ > 0$. We make the following additional assumptions on $\htp{\lambda}$ and $\mu$. 
\begin{assum}\label{assum:heavy_traffic} All of the following hold for arrival rates $\htp{\lambda}$ given by \eqref{eq:customer_arrival_rate} and service rates $\mu$:
\begin{enumerate}[nosep]
    \item[\rm (i)] $|\Lambda| = |\mu|$,
    \item[\rm (ii)] $|\gamma| > 0$,
    \item[\rm (iii)] $\gamma_i < 0$ for all $i \in [n]$ such that $\Lambda_i = 0$.
\end{enumerate}
\end{assum}
Parts (i) and (ii) ensure that we are approaching heavy traffic from below. Part (iii) is implied by $\htp{\lambda_i} > 0$ for all $i \in [n]$ and $0 < \epsilon < \epsilon_+$, but we include it in \cref{assum:heavy_traffic} for clarity. Note that for $i \in [n]$ such that $\Lambda_i > 0$, we allow $\gamma_i$ to be positive, negative, or zero. 

This is more general than the scaling used in \cite{Afecheetal2019}, where the authors assume that $\gamma = \Lambda$. Additionally, \cite{Afecheetal2019} requires that $\Lambda_i > 0$ for all $i \in [n]$. We relax that assumption here, as it is useful to allow for no arrivals to particular customer classes when considering strategic customer behaviour. 

We are only interested in studying systems which produce stable outcomes. This leads us to restrict our attention to a set of {\em admissible} menus. 
\begin{dfn}{\rm (Admissible Menus)}\label{def:admissible}
For a given menu $M$, arrival rates $\htp{\lambda}$, and service rates $\mu$, define for any subset of servers $\mathscr{S} \subseteq [m]$ and $\epsilon > 0$
\begin{equation}
\htp{\Delta}_{\mathscr{S}}(M) := \sum_{j \in \mathscr{S}} \mu_j - \sum_{i \in U_{\mathscr{S}}(M)} \htp{\lambda}_i.
\end{equation}

A menu M is {\rm admissible} under arrival rates  $\htp{\lambda}$ and service rates $\mu$ if 
\begin{equation*}
    \htp{\Delta}_{\mathscr{S}}(M) = \Omega(\epsilon) \qquad \mbox{for all } \mathscr{S} \subseteq [m].
\end{equation*}
\end{dfn}
In words, this ensures that the menu $M$ and arrival rates $\htp{\lambda}$ admit a steady state under a FCFS-ALIS service discipline, and that the slack in the system is converging slowly enough so that the average delays of the different customer classes converge when scaled by $\epsilon$. 

We let $\cal{M}(\htp{\lambda}, \mu)$ denote the set of all menus $M$ that are admissible for arrival rates $\htp{\lambda}$ and service rates $\mu$. The set $\cal{M}(\htp{\lambda}, \mu)$ will be non-empty for all pairs $(\htp{\lambda}, \mu)$ satisfying \cref{assum:heavy_traffic}. To see this, observe that the complete menu $M$ such that $m_{ij} = 1$ for all $i \in [n]$ and $j \in [m]$ will be admissible for all $(\htp{\lambda}, \mu)$ satisfying \cref{assum:heavy_traffic}. The complete menu will operate like a single queue with arrival rates $|\htp{\lambda}|$ that is served by all servers. 

\section{Mean Waiting Times in Heavy Traffic}\label{sec:heavytrafficwaits}
We are interested in being able to calculate the mean waiting times of the different service classes. Because we are looking at a conventional heavy traffic setting, the waiting times themselves will grow out of bound as $\epsilon \downarrow 0$. We will instead look at the scaled mean waiting time 
\begin{equation}
    \hts{W_i} = \epsilon \cdot \htp{W_i},
\end{equation}
which will remain bounded in heavy traffic. 

In what follows we show how to find the limiting expected waiting times by building upon and extending the methods and results
in \cite{Afecheetal2019}.

\subsection{Feasible flows and CRP components}
We begin by identifying the feasible flows of customers between customer classes and servers. For arrival rates $\htp{\lambda}$ and service rates $\mu$ satisfying \cref{assum:heavy_traffic}, and an admissible menu $M \in \cal{M}(\htp{\lambda}, \mu)$, for $0 \le \epsilon < \epsilon_0$ we define the set of feasible flows as
\begin{align}
    {\cal F}(\epsilon, \htp{\lambda}, M) := \Big\{f=[f_{ij}] \geq 0:  & \sum_{i \in [n]} f_{ij} \le \mu_j,~~\forall j \in {[m]}; \nonumber \\ 
    & \sum_{j \in [m]} f_{ij} = \htp{\lambda_i},~~\forall i \in [n];~~ f_{ij}= 0, ~~\forall (i,j) : m_{ij}=0 \Big\},
\end{align}
where $\epsilon_0 \in \RR$ is such that $\htp{\lambda} > 0$ for all $0 < \epsilon < 0$. We know from the admissibility of $M$ that such an $\epsilon_0$ exists, and that $\mathcal{F}(\epsilon, \htp{\lambda}, M)$ is non-empty for all $0 < \epsilon < \epsilon_0$. The following lemma shows that $ {\cal F}(0, \htp{\lambda} M)$ is also non-empty. The proof relies on $ {\cal F}(\epsilon, \htp{\lambda}, M)$ being a subset of a compact set ${\cal F}_{\max}(\htp{\lambda})$ for $0 \le \epsilon < \epsilon_0$. 

\begin{lem}\label{lem:flows}
For a given $\htp{\lambda}$ and $\mu$ satisfying 
\cref{assum:heavy_traffic}, and $M \in \cal{M}(\htp{\lambda}, \mu)$, the set ${\cal F}(0, \Lambda, M)$ is non-empty. Furthermore, every sequence of flows $\htp{f}$ such that $\htp{f} \in {\cal F}(\epsilon, \htp{\lambda}, M)$ has a sub-sequence that converges to some $\tilde{f} \in {\cal F}(0, \Lambda,  M)$.  
\end{lem}
{\sc Proof}: See Appendix \ref{sec:proof_W}. $\Box$

As this lemma suggests, the set ${\cal F}(0, \Lambda, M)$ contains information about what sort of flows it is possible to observe in heavy traffic. We will use the set of feasible limiting flows to determine which servers have a positive probability of serving which service classes in the limit. To do this, we will first define the \textit{residual matching} of the menu $M$. 

\begin{dfn}\label{def:residualmatching} {\rm (Residual Matching)} For a given $(\htp{\lambda},\mu,M)$ such that $\htp{\lambda}$ and $\mu$ satisfy 
\cref{assum:heavy_traffic} and $M \in \cal{M}(\htp{\lambda}, \mu)$ we define the {\em residual matching} $\breve{M}$, where $\breve{M}=[\breve{m}_{ij}]$ satisfies $\breve{m}_{ij}=1$ if and only if there exists flows $\tilde{f} \in {\cal F}(0,\Lambda, M)$ such that $\tilde{f}_{ij} > 0$.
\end{dfn} 
Intuitively, for a service class $i$ and server $j$ with $m_{ij} = 1$ but  $\breve{m}_{ij}=0$, the flow of customers from service class $i$ to server $j$ must vanish in the heavy-traffic limit. 
\cite{Afecheetal2019} provide an algorithm for finding the residual matching. However, for small, simple systems the residual matching can be found by inspection. To see this, consider again the simple example in \cref{fig:Ex1}, specifying the service rates to be $\mu = [2, 1, 2, 1]$. We will consider two example vectors of arrival rates, $\Lambda_a = [2, 1, 1, 2]$ and $\Lambda_b = [2, 1, 0, 3]$. In each case, there is only one set of feasible flows in  ${\cal F}(0, \Lambda_a, M)$ and  ${\cal F}(0, \Lambda_b, M)$, given by
\begin{equation}
f^a_{ij}=\left[\begin{array}{cccc} 2 & 0 & 0 & 0 \\ 0 & 1 & 0 & 0 \\ 0 & 0 & 1 & 0 \\ 0 & 0 & 1 & 1\end{array}\right] \quad \mbox{and} \quad f^b_{ij}=\left[\begin{array}{cccc} 2 & 0 & 0 & 0 \\ 0 & 1 & 0 & 0 \\ 0 & 0 & 0 & 0 \\ 0 & 0 & 2 & 1\end{array}\right].
\end{equation}

\begin{figure}[hbt]
    \begin{center}
         \begin{subfigure}[b]{0.45\textwidth}
             \centering
             \includegraphics[width=7cm]{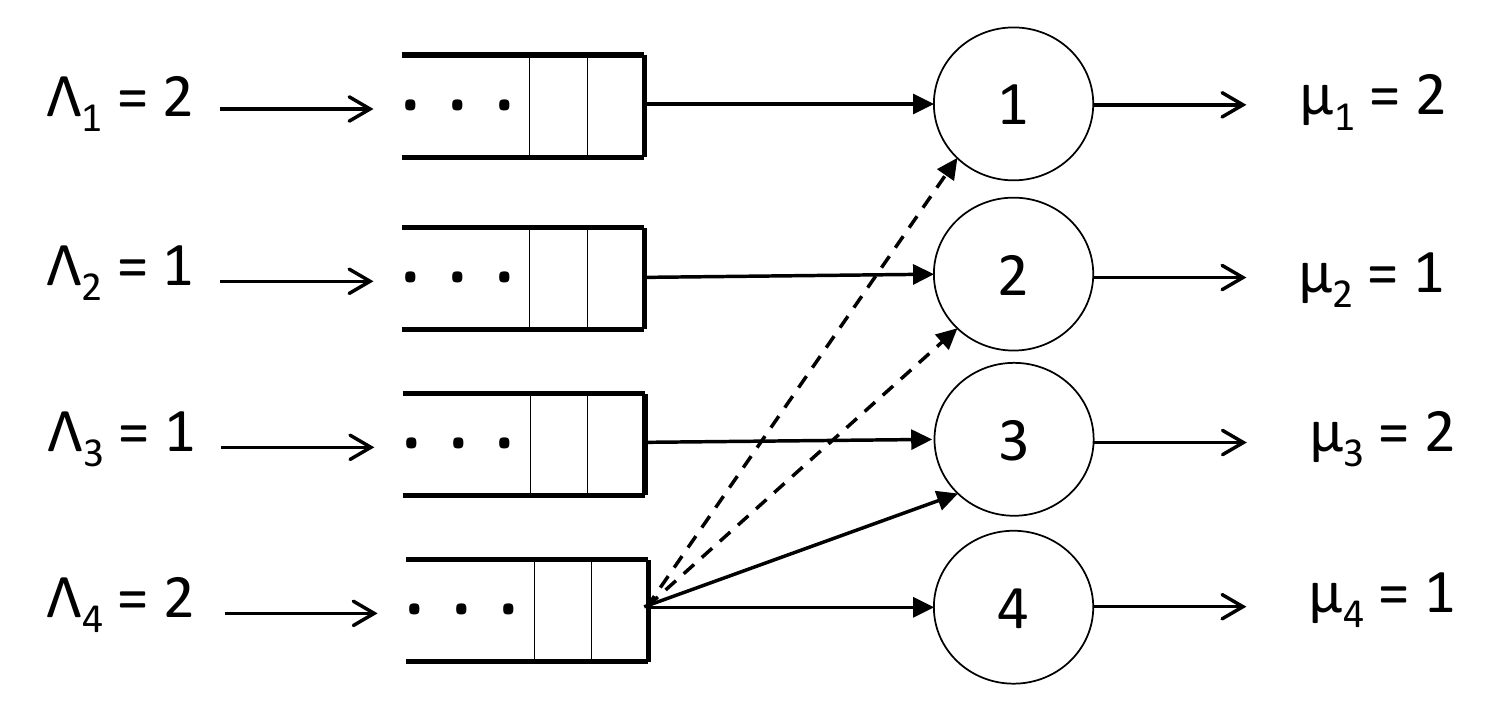}
             \caption{Residual matching (a)}
             \label{subfig:residualmatching_a}
         \end{subfigure}
         \hfill
         \begin{subfigure}[b]{0.45\textwidth}
             \centering
             \includegraphics[width=7cm]{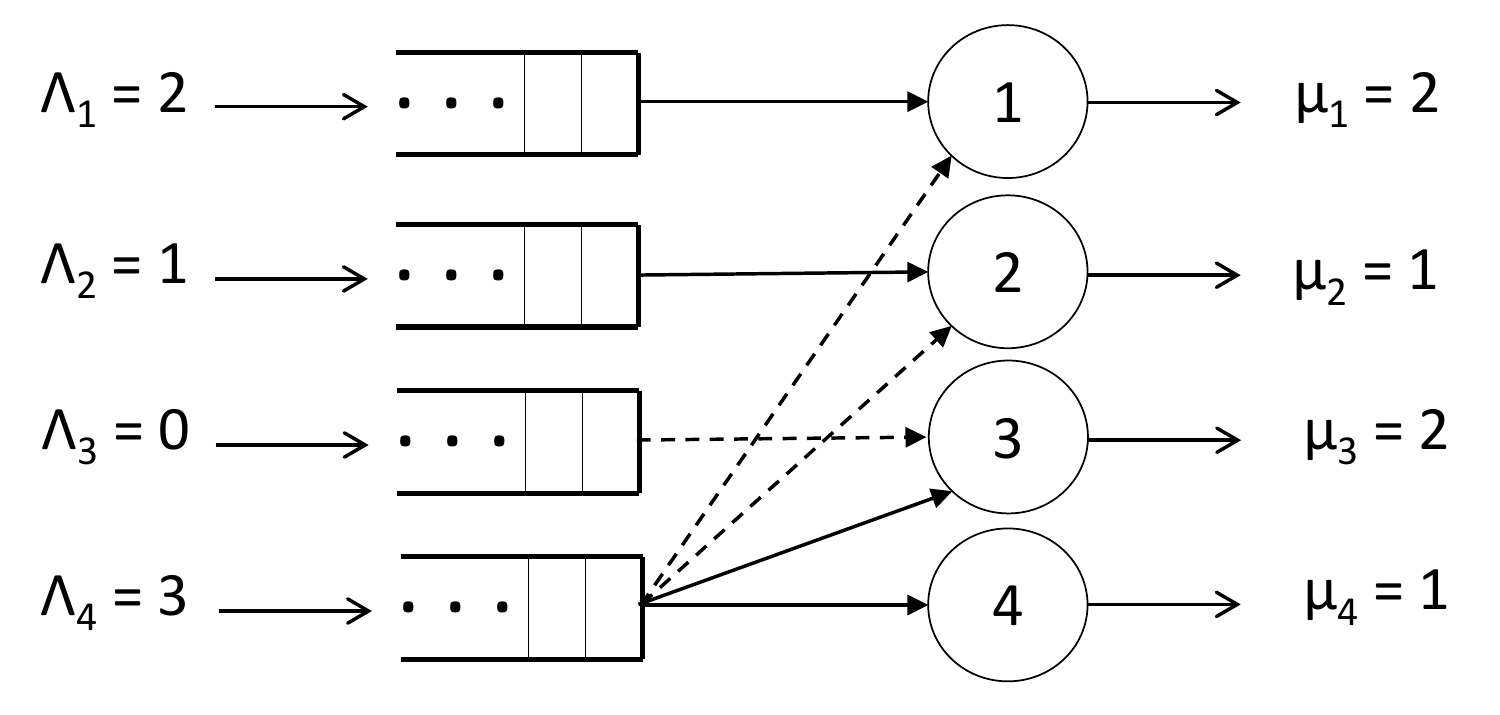}
             \caption{Residual matching (b)}
             \label{subfig:residualmatching_b}
         \end{subfigure}
    \caption[Model Description]{\sf \footnotesize Examples of residual matchings.}
     \label{fig:residual_matching}
    \end{center}
\end{figure}

In example (a), the arcs in the compatibility network with $m_{ij} = 1$ and $\breve{m}_{ij} = 0$ are (4,1) and (4,2). While service class 4 is compatible with servers 1 and 2, there will be zero flow between class 4 and servers 1 and 2 in the limit. All the service capacity of servers 1 and 2 will be allocated to serving classes 1 and 2. We can see this visually in panel (a) in \cref{fig:residual_matching}, where the arcs with $m_{ij} = 1$ and $\breve{m}_{ij} = 1$ are represented with solid lines, and the arcs with $m_{ij} = 1$ and $\breve{m}_{ij} = 0$ are represented with dashed lines. Example (b) is similar, but we now additionally have arc (3,3) with $m_{33} = 1$ and $\breve{m}_{33} = 0$. In panel (b) of \cref{fig:residual_matching} we can see that class 3 only has one dashed arc connecting it to any servers, representing that no servers are allocating any capacity to class 3 in the limit, even though class 3 is compatible with server 3. 

Knowing the residual matching allows us to decompose the initial bipartite matching system into a partition of independent components, which \cite{Afecheetal2019} refer to as {\em complete resource pooling} (CRP) components. 

\begin{dfn}\label{def:CRP}{\rm (CRP Component)}  For a given $(\htp{\lambda},\mu,M)$ such that $\htp{\lambda}$ and $\mu$ satisfy 
\cref{assum:heavy_traffic} and $M \in \cal{M}(\htp{\lambda}, \mu)$, let the induced residual matching be denoted $\breve{M}$. We say that the subset  \linebreak $\mathbb{C}=(\mathcal{C}, \mathcal{S}) \in 2^{[n]}\times 2^{[m]}$ of customer classes and servers forms a \emph{CRP component} if for any pair of nodes $k_1, k_2 \in \mathcal{C} \cup \mathcal{S}$ there exists a path between $k_1$ and $k_2$ in $\breve{M}$, and $\mathbb{C}$ is maximal in the sense that the condition is violated for any strict superset of $\mathbb{C}$.
\end{dfn}

We let $\{\mathbb{C}_1,\mathbb{C}_2,\dots,\mathbb{C}_K\}$ denote the collection of CRP components induced by the residual matching $\breve{M}$, where $K$ is the number of components. Each $\mathbb{C}_k = (\mathcal{C}_k, \mathcal{S}_k)$ is defined by the subset of customer classes $\mathcal{C}_k$ and the subset of servers $\mathcal{S}_k$ that belong to $\mathbb{C}_k$. Since we allow for service classes with no arrivals, that is $\Lambda_i = 0$, some CRP components will have an empty server set. Each service class with $\Lambda_i = 0$ forms a separate CRP component with an empty server set. We denote the subset of such CRP components by $\Ical_0$:
\begin{align}
        \Ical_{0} &= \{ k : \Lambda_k = 0\}.
\end{align}
We let $K' := K - |\Ical_0|$ be the number of CRP components with non-empty sets of servers, and will assume that the CRP components are indexed so that the components in $[K] \setminus \Ical_0$ have indices $1,2, \ldots, K'$. We will use $k(i)$ and $k(j)$ to denote the component that service class $i$ or server $j$ is part of, where the use should be clear from context.

To make these ideas more concrete, let us return to our examples in \cref{fig:residual_matching}. In example (a), service class 1 and server 1 make up a CRP component, as they are not connected to any other service classes or servers with solid arcs. Similarly, service class 2 and server 2 make up a CRP component. We can see a path between classes 3 and 4 through server 3, so these classes along with servers 3 and 4 make up a single CRP component. This means the CRP components for example (a) can be written as $\mathbb{C}_1 = (\mathcal{C}_1, \mathcal{S}_2 = (\{1\}, \{1\})$, $\mathbb{C}_2 = (\mathcal{C}_2, \mathcal{S}_2 = (\{2\}, \{2\})$, and $\mathbb{C}_3 = (\mathcal{C}_3, \mathcal{S}_3 = (\{3, 4\}, \{3, 4\})$. Example (b) is similar, the difference being that now service class 3 is not connected to any server or service class with a solid arc, and therefore is in a CRP component by itself with an empty server set, i.e. ${\cal I}_0=\{3\}$. So the CRP components for example (b) are $\mathbb{C}_1 = (\mathcal{C}_1, \mathcal{S}_2 = (\{1\}, \{1\})$, $\mathbb{C}_2 = (\mathcal{C}_2, \mathcal{S}_2 = (\{2\}, \{2\})$, $\mathbb{C}_3 = (\mathcal{C}_3, \mathcal{S}_3 = (\{4\}, \{3, 4\})$, and $\mathbb{C}_4 = (\mathcal{C}_4, \mathcal{S}_4 = (\{3\}, \{ \emptyset \})$. 

Abusing notation, we denote the aggregate arrival and service rates for the CRP components under $\htp{\lambda}$ as: 
\begin{align}
\label{eq:CRP_mu_lambda}
\forall k \in [K] & \ : \ \htp{\widetilde{\lambda}_k} = \sum_{i \in \Ccal_k} \htp{\lambda_i} =: \widetilde{\Lambda}_k - \epsilon \widetilde{\gamma}_k  + o(\epsilon),\quad \mbox{and} \quad  \widetilde{\mu}_k = \sum_{j \in \Scal_k} \mu_j , 
\end{align}
where $\widetilde{\Lambda}_k = \sum_{i \in \Ccal_k} \Lambda_i$ and $\widetilde{\gamma}_k = \sum_{i \in \Ccal_k} \gamma_i$. We will later show that each CRP component must satisfy $\widetilde{\Lambda}_k = \widetilde{\mu}_k$ so that the slack between demand and capacity within a CRP component in heavy-traffic goes to zero with $\epsilon$. 
While each CRP component is critically loaded, the ``well-connectedness'' within a CRP component allows shifting load from one service class to another on short time scales. In particular, we will show in \cref{thm:CRPdelay_append} that under an FCFS-ALIS policy, waiting times are balanced in such a way that service classes that belong to the same CRP component have the same limiting scaled mean waiting time in the heavy traffic limit. 

\subsection{Directed Acyclic Graph of CRP components}
The menu $M$ and the residual matching $\breve{M}$ uniquely induce a directed acyclic graph (DAG) on the collection of CRP components defined in the previous step. This is useful as the DAG defines a precedence relation among service classes: since component $k_1$ has a directed arc to component $k_2$, there is a service class in $k_1$ that can be served by a server in $k_2$. This means $k_1$ can ``off-load'' its customers to the servers of component $k_2$, and so the instantaneous waiting time in component $k_1$ cannot exceed that in component $k_2$ under FCFS-ALIS. This intuition is made precise in the proof of \cref{thm:CRPdelay_append}.

The following is a formal statement of how the DAG is induced. 

\begin{dfn}\label{def:DAG}{\rm (DAG)} Given the menu $M=[m_{ij}]$, and the CRP components $\{\mathbb{C}_k=(\mathcal{C}_k, \mathcal{S}_k) \colon k=1,\dots,K\}$ induced by the residual matching $\breve{M}$, we define  $\Dcal = ([K], \Acal)$ associated to $M$ as the directed acyclic graph whose nodes correspond to the CRP components, and there is a directed arc $(k_1, k_2) \in \Acal$ from component $\mathbb{C}_{k_1}$ to component $\mathbb{C}_{k_2}$ if and only if there exists a customer class $i \in  \mathcal{C}_{k_1}$ and a server $j \in \mathcal{S}_{k_2}$ such that $m_{ij}=1$. We use the notation $k_1 \stackrel{\Dcal}{\rightsquigarrow} k_2$ to denote that there is a directed path $k_1$ to $k_2$ in the DAG $\Dcal$.
\end{dfn}

\cite[Lemma 2]{Afecheetal2019} formally proves that the directed graph defined above is in fact acyclic. 

Returning to our examples in \cref{fig:residual_matching}, the DAGs are given below. 
\begin{figure}[hbt]
    \begin{center}
         \begin{subfigure}[b]{0.45\textwidth}
             \centering             \includegraphics[width=6cm]{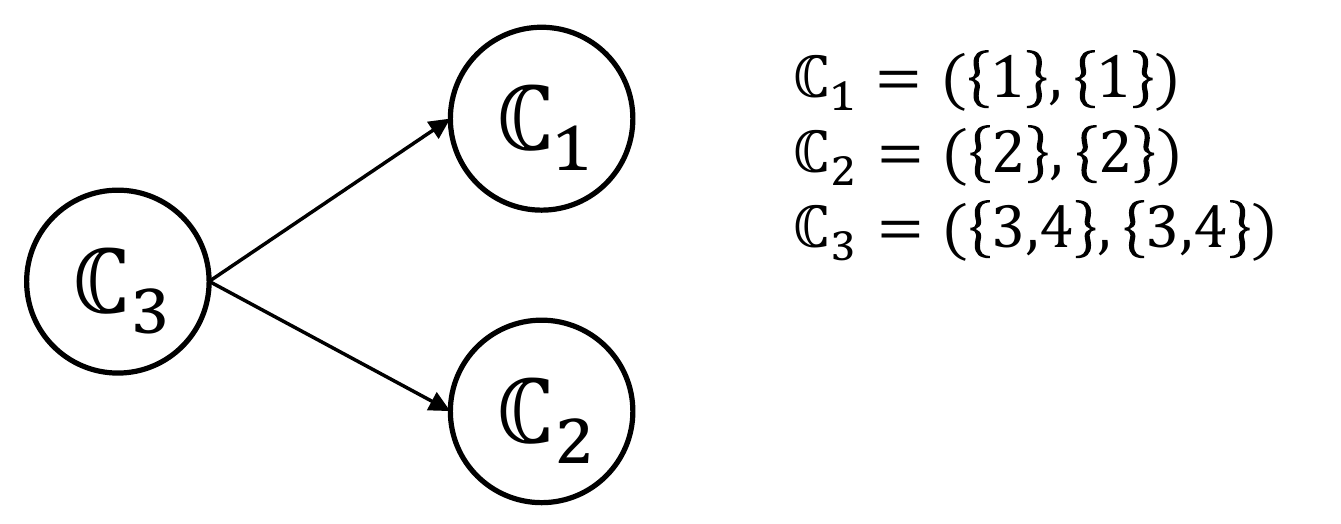}
             \caption{DAG (a)}
             \label{subfig:DAG_a}
         \end{subfigure}
         \hfill
         \begin{subfigure}[b]{0.45\textwidth}
             \centering
             \includegraphics[width=6cm]{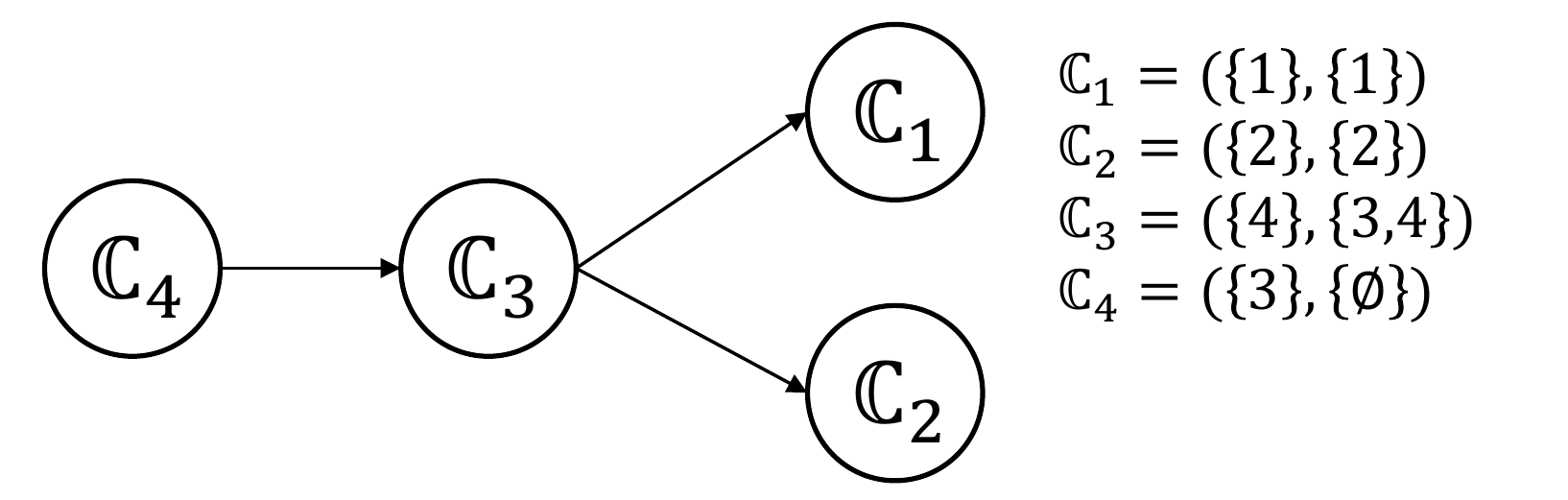}
             \caption{DAG (b)}
             \label{subfig:DAG_b}
         \end{subfigure}
    \caption[Model Description]{\sf \footnotesize Examples of DAGs.}
     \label{fig:DAGs}
    \end{center}
\end{figure}
In both cases, service class 4 can be served by servers 1 and 2 in the original menu, i.e. $m_{41} = m_{42} = 1$, and so there are directed arcs from $\mathbb{C}_3$ to $\mathbb{C}_1$ and $\mathbb{C}_2$. In example (b), $\mathbb{C}_4$ contains service class 3 but no servers, since service class 3 has an arrival rate of 0. Therefore $\mathbb{C}_4$ has a directed arc to $\mathbb{C}_3$, as this is the CRP component containing the server that customer class 3 is compatible with. 

As we mentioned earlier, our computations for the heavy-traffic waiting times build on the work of \cite{AdanWeiss2014}. The crucial component of their analysis is a state-space representation for the FCFS-ALIS matching model which involves ranking the busy servers in order of the waiting time of the customers they are serving. As was proved in \cite{Afecheetal2019} for the less general scaling, in heavy-traffic this entails restricting attention to only certain permutations of the CRP components which have asymptotically non-zero steady-state probability. We show in  \cref{prop:steadystate} below that this also holds for our more general scaling. The topological orders of the DAG $\Dcal$ are precisely these permutations. The definition we give next differs slightly from \cite{Afecheetal2019} due to the potential presence of CRP components with $\widetilde{\Lambda}_k=0$.

\begin{dfn}\label{def:topologicalorders}{\em (Topological Orders on CRP Components)} Let $\{\mathbb{C}_1, \mathbb{C}_2,\dots,\mathbb{C}_{K'}\}$ be the CRP components with $\widetilde{\Lambda}_k > 0$. Given the DAG $\Dcal = ([K],\Acal)$, we say that a permutation $\sigma=(\sigma(1),\sigma(2),\dots,\sigma(K'))$ of $[K']$ induces a {\rm topological order} $(\mathbb{C}_{\sigma(1)}, \mathbb{C}_{\sigma(2)},\dots,\mathbb{C}_{\sigma(K')})$ of these CRP components if for every pair $(k_1,k_2) \in [K']$ such that $k_1 \stackrel{\Dcal}{\rightsquigarrow} k_2$, we have $\sigma^{-1}(k_2) < \sigma^{-1}(k_1)$. In other words, sink components of $\Dcal$ precede source components. We let ${\cal T}(\Dcal, K')$ denote the set of all permutations $\sigma$ of $[K']$ that induce a topological order on components $\{\mathbb{C}_1,\ldots,\mathbb{C}_{K'}\}$.

Further, for each $\sigma \in \Tcal(\Dcal, K')$, we partition the CRP components $[K]$ by associating a subset for each $k \in [K']$ as follows:
\begin{align}
\label{eq:comps_fn}
    \comps(\sigma, k) := \{\sigma(k)\} \union \{ \kappa \in \Ical_0 : k = \max \{ k' \in [K'] : \  \kappa \stackrel{\Dcal}{\rightsquigarrow} \sigma(k') \} \}.
\end{align}
The interpretation of this is that for each index $k \in [K']$, we associate the CRP component corresponding to $\sigma(k)$ as well as all CRP components $\kappa$ with $\widetilde{\Lambda}_\kappa = 0$ (i.e., server-less components) for which the component $\sigma(k)$ is the last component in the topological order $\sigma$ that is reachable from $\kappa$ via a directed path. 
\end{dfn}
We will use the shorthand $\comps^{-1}(\sigma, k)$ to denote the index $\kappa \in [K']$ such that $k \in \comps(\sigma, \kappa)$.

To highlight the difference with \cite{Afecheetal2019}, under the heavy-traffic regime considered in \cite{Afecheetal2019} all CRP components have a non-empty server set $\Scal_j$. In contrast, in our model, we have customer classes that are in CRP components by themselves. These CRP components are special in that they have no incoming arc in the DAG $\Dcal$, and can only have a directed arc to CRP components with non-empty server sets. The topological orders $\Tcal(\Dcal,K')$ can thus be thought of as preprocessing $\Dcal$ to remove the server-less CRP components $\{\mathbb{C}_{K'+1},\ldots, \mathbb{C}_K\}$ which are ``hanging off'' $\Dcal$, and finding topological orders on the remaining components. Since the topological order has sink components of $\Dcal$ preceding source components, and as we mentioned earlier, the DAG defines a precedence relation among service classes, we can then interpret $\comps^{-1}(\sigma, k)$ as associating each server-less CRP component with the CRP component that is reachable from it that has the shortest steady-state wait. 

Returning to our examples in \cref{fig:DAGs}, both example (a) and example (b) have the same set of CRP components with positive limiting arrival rates, the set $\{\mathbb{C}_1, \mathbb{C}_2, \mathbb{C}_3  \}$. Both examples also have the same connectivity with these components. $\mathbb{C}_3$ has directed arcs to $\mathbb{C}_1$ and $\mathbb{C}_2$, but there are no arcs between $\mathbb{C}_2$ and $\mathbb{C}_2$. Hence in any topological orders on these CRP components, we know that $\mathbb{C}_1$ and $\mathbb{C}_2$ come before $\mathbb{C}_3$, but $\mathbb{C}_1$ can come either before or after $\mathbb{C}_2$. Thus the possible permutations are $\sigma_1 = (1, 2, 3)$ and $\sigma_2 = (2, 1, 3)$, and the  associated topological orders are $(\mathbb{C}_1, \mathbb{C}_2, \mathbb{C}_3)$ and $(\mathbb{C}_2, \mathbb{C}_1, \mathbb{C}_3)$. As example (a) has no CRP components with limiting arrival rates of 0, for each $\sigma$ and each $k$, $\comps(\sigma, k)$ is simply the set containing the index of the CRP component at position $k$ of the topological order $\sigma$. In example (b), $\mathbb{C}_4$ has $\tilde{\lambda}_4 = 0$, so for each topological order $\sigma$, we need to determine for which $k$ we have $4 \in \comps(\sigma, k)$. The only directed arc from $\mathbb{C}_4$ to any other CRP component is to $\mathbb{C}_3$. Hence for each $\sigma$, we have that $4 \in \comps(\sigma, k)$ if and only if $3 \in \comps(\sigma, 4)$. Since $\mathbb{C}_3$ is the last element of the topological order for both permutations $\sigma_a$ and $\sigma_b$, we have that $\comps(\sigma_a, 3) = \comps(\sigma_b, 3) = \{3, 4\}$. 

\subsection{Calculating waiting times}

Let $\Tcal(\Dcal, K') = (\sigma_1, \ldots, \sigma_T)$ be the collection of topological orders on $\{\mathbb{C}_1, \ldots, \mathbb{C}_{K'}\}$ (the components with $\widetilde{\Lambda}_k > 0$). For a topological order $\sigma_t \in \Tcal(\Dcal, K')$ with the associated function $\comps(\sigma_t, \cdot)$ defined in \eqref{eq:comps_fn}, we define the unnormalized probability of being in a state associated with the topological order $\sigma_t$ as:
\begin{align}
\label{eq:Q_def}
    \mathbb{Q}(\sigma_t) &= \prod_{\kappa \in [K']} \frac{1}{\sum_{\ell=1}^\kappa  \widetilde{\gamma}_{\comps(\sigma_t, \ell)}},
\end{align}
where we use the shorthand 
\[ \widetilde{\gamma}_{\comps(\sigma, \ell)} = \sum_{\kappa \in \comps(\sigma, \ell)} \widetilde{\gamma}_\kappa .\]
For a permutation $\sigma_t \in \Tcal(\Dcal, K')$, for any CRP component $\mathbb{C}_k$, we define the waiting time conditioned on the topological order $\sigma_t$ as:
\begin{align}
\label{eq:w_sigma_k}
    w_{\sigma_t, k} &= \sum_{ \kappa =  \comps^{-1}(\sigma_t, k)  }^{K'} \frac{1}{\sum_{\ell=1}^\kappa \widetilde{\gamma}_{\comps(\sigma_t, \ell)}} .
\end{align}
The following \cref{lem:slack_positive} proves that the expressions above are well-defined.
\begin{lem}
\label{lem:slack_positive}
For $\htp{\lambda}$ and $\mu$ satisfying \cref{assum:heavy_traffic}, and for some $M \in \cal{M}(\htp{\lambda}, \mu)$ for all permutations $\sigma_t \in \Tcal(\Dcal, K')$ of CRP components $\{\mathbb{C}_1,\ldots,\mathbb{C}_{K'}\}$ and for all $\kappa \in [K']$,
\[  \sum_{\ell = 1}^\kappa \widetilde{\gamma}_{\comps(\sigma_t, \ell)} > 0. \]
    \end{lem}
{\sc Proof}: See Appendix \ref{sec:proof_W}. $\Box$
    
With the expressions for the unnormalized probabilities and conditional waiting times of topological orders in place, we are ready to state our main theorem regarding the mean scaled steady-state waiting times of different service classes. 

\begin{thm}\label{thm:CRPdelay_append} For a given $(\htp{\lambda},\mu,M)$ such that $\htp{\lambda}$ and $\mu$ satisfy 
\cref{assum:heavy_traffic}, and an admissible menu $M \in \cal{M}(\htp{\lambda}, \mu)$, let $\breve{M}$ be the residual matching and  $\{\mathbb{C}_1,\dots,\mathbb{C}_{K'}, \mathbb{C}_{K'+1}, \dots, \mathbb{C}_{K} \}$ be the collection of CRP components induced by $\breve{M}$. Then, customer classes that belong to the same CRP component experience the same scaled steady-state mean waiting time in heavy traffic. Furthermore, the scaled steady-state mean waiting time of CRP component  $\mathbb{C}_k$ is equal to
\begin{align} 
\label{eq:WCRP_def}
\widehat{W}_{\mathbb{C}_k} = \sum_{t=1}^{T(M)} \left({\mathbb{Q}(\sigma_t) \over \mathbb{Q}(\sigma_1)+\mathbb{Q}(\sigma_2)+\cdots+\mathbb{Q}(\sigma_{T(M)})}\right)\, w_{\sigma_t, k}.
\end{align}
\end{thm}

The proof of \cref{thm:CRPdelay_append} can be found in \cref{sec:waits_proof}.

\section{Matching Probabilities in Heavy Traffic}\label{sec:matchingprob}
Another performance metric of interest is the matching probabilities, that is, for each customer class $i$ and server $j$, the probability that a customer who joins class $i$ is served by server $j$. For any menu $M$ that is admissible with arrival rates $\htp{\lambda}$ and service rates $\mu$, we let $\htp{p}(M, \htp{\lambda}, \mu)$ be the matrix of matching probabilities, so $\htp{p}_{ij}(M, \htp{\lambda}, \mu)$ is the steady state probability with which a customer who joins class $i \in [n]$ is served by server $j\in [m]$.  While exact matching probabilities are difficult to calculate, and remain difficult to calculate even in heavy traffic, we are able to provide two results regarding how matching rate calculations simplify as we move to heavy traffic.  

Before stating our results, it will be useful to describe the combinations of limiting arrival rates $\Lambda$, service rates $\mu$, and menus $M$ such that there is some sequence $\htp{\lambda}$ converging to $\Lambda$ that makes $M$ admissible. The following proposition will help us understand these combinations. 
\begin{prop}\label{prop:limiting_admissibility} Take any sequence of arrival rates $\htp{\lambda}$ and service rates $\mu$ such that , $\htp{\lambda}$ and $\mu$ satisfy 
\cref{assum:heavy_traffic}, and let $M$ be such that $M \in \cal{M}(\htp{\lambda}\mu)$. Let $\Lambda = \lim_{\epsilon \to 0}\htp{\lambda}$. Then $M$ is admissible with service rates $\mu$ and arrival rates 
\begin{equation*}
    \htp{\lambda} = \Lambda - \epsilon\Lambda, \quad \text{for } \epsilon > 0. 
\end{equation*}
Furthermore, if $M$ is admissible with $\htp{\lambda} = \Lambda - \epsilon\Lambda$ and $\mu$, then the menu $\breve{M}$ given by the residual matching of $M$ is also admissible with $\htp{\lambda} = \Lambda - \epsilon\Lambda$ and $\mu$. 
\end{prop}
{\sc Proof}: See Appendix \ref{sec:proof_matchingprob_statement}. $\Box$

This lets us talk about menus that are admissible for limiting arrival rates $\Lambda$ and service rates $\mu$. We will define the set $\mathcal{M}^+(\Lambda, \mu)$ to be the set of all menus $M$ such that $M$ is admissible for arrival rates $\htp{\lambda} = \Lambda(1 - \epsilon)$ and service rates $\mu$. This provides us with a more convenient way to express our results regarding matching probabilities, the first of which is stated formally in \cref{thm:matching_independent_gamma}. This tells us that while the limiting expected delays depend on the particular sequence of arrival rates $\htp{\lambda}$, and in particular depend on the slacks $\gamma$, the matching probabilities depend only on the limiting arrival rates.

\begin{thm}\label{thm:matching_independent_gamma} Take any limiting arrival rates $\Lambda$ and service rates $\mu$ such that $|\Lambda| = |\mu|$. Consider any menu $M \in \mathcal{M}^+(\Lambda, \mu)$. Take any two sequences of arrival rates $\htp{\lambda}_a$ and $\htp{\lambda}_b$ such that $\lim_{\eps \to 0}\htp{\lambda}_a = \lim_{\eps \to 0}\htp{\lambda}_b = \Lambda$, both sequences satisfy \cref{assum:heavy_traffic} with $\mu$, and $M$ is admissible for both sequences of arrival rates with $\mu$. Then $\lim_{\epsilon \to 0}\htp{p}_{ij}(M, \htp{\lambda}_a, \mu) = \lim_{\epsilon \to 0}\htp{p}_{ij}(M, \htp{\lambda}_b, \mu)$ for all $i \in [n]$ and $j \in [m]$. 
\end{thm}
\cref{thm:matching_independent_gamma} and \cref{thm:independent_matching_rates} can be found in \cref{sec:proof_matching}. 

\cref{thm:matching_independent_gamma} lets us talk about the matching probabilities of a menu $M$ just in terms of the limiting arrival rates $\Lambda$ and service rates $\mu$. In light of this, for the rest of this paper we will refer to matching probabilities in terms of the limiting arrival rates, that is, we will write $\htp{p}_{ij}(M, \Lambda, \mu)$. 
The proof of 

The second result we have relating to matching probabilities, stated formally in \cref{thm:independent_matching_rates}, tells us that matching probabilities within a CRP component are independent of all other CRP components. 
\begin{cor}\label{thm:independent_matching_rates} Take any limiting arrival rates $\Lambda$ and service rates $\mu$ such that $|\Lambda| = |\mu|$, and take any $M \in \mathcal{M}^+(\Lambda, \mu)$. Let $\breve{M}$ be the residual matching, and let $\{\mathbb{C}_1,\dots,\mathbb{C}_{K'}, \mathbb{C}_{K'+1}, \dots, \mathbb{C}_{K} \}$ be the collection of CRP components induced by $\breve{M}$. Then for any customer class $i \in \Ccal_k$ and server $j \in \Scal_k$, 
\begin{equation*}
\lim_{\epsilon \to 0}\htp{p}_{ij}(M, \Lambda, \mu) = \lim_{\epsilon \to 0}\htp{p}_{ij}(\breve{M}, \Lambda, \mu).
\end{equation*}
\end{cor}

\cref{thm:independent_matching_rates} implies that when calculating the matching rates, we can look at each CRP component individually. Additionally, it tells us that the DAG structure does not affect the matching probabilities. We will see in \cref{sec:discussion} that two menus $M$ and $M'$ with the same residual matching $\breve{M}$ can have significantly different expected waiting times in heavy-traffic if the two menus induce different DAGs. \cref{thm:independent_matching_rates} tells us that despite this, the limiting matching probabilities of menus $M$ and $M'$ are the same.

\section{Discussion}\label{sec:discussion}
Before getting into the proofs of our main results, we discuss some of their implications, while highlighting the differences between the behaviours of our model and the model in \cite{Afecheetal2019}. We also explore some simple questions regarding the design of menus of service classes. 

\subsection{Implementable outcomes}
Our motivation for the heavy-traffic scaling used in this paper is that it allows for a wider range of outcomes than the proportional scaling used in \cite{Afecheetal2019}. The following definition will help formalise what we mean by this. 

\begin{dfn}\label{def:WImplement}{\rm (Implementable Waiting Times)} Take limiting arrival rates $\Lambda$, service rates $\mu$, and a menu $M$ such that a collection of CRP components $\mathbb{C}=\{\mathbb{C}_1,\mathbb{C}_2,\dots,\mathbb{C}_K\}$ is induced. We say a vector of limiting scaled waiting times $W=(W_1,W_2,\dots,W_K)$ is {\rm implementable} if there exists $\gamma \in \RR^n$ such that the menu $M$ is admissible for the pair $(\htp{\lambda}, \mu)$ where 
\begin{equation*}
    \htp{\lambda}_i = \Lambda_i - \epsilon\gamma_i + o(\epsilon), \quad \text{for all } i \in [n],
\end{equation*}
and the resulting limiting waiting times $\widehat{W}_{\mathbb{C}_k}$ given by \eqref{eq:WCRP_def} are equal to $W_k$ for all $k\in [K]$.
\end{dfn}

If we only look at the scaling in \cite{Afecheetal2019}, in which $\gamma = \Lambda$, then each combination of limiting arrival rates $\Lambda$, service rates $\mu$, and menu $M$ can produce one specific vector of waiting times. By allowing $\gamma$ to change, we increase the set of implementable outcomes. 

As we alluded to in \cref{sec:heavytrafficwaits}, the DAG provides information about which vectors of waiting times are implementable. The following statement, which is a corollary of \cref{thm:CRPdelay_append}, formalises this idea. 

\begin{cor}\label{cor:implementable}
    If $W \in \RR^K_+$ is implementable, then $W$ is consistent with some topological order $\sigma \in \mathcal{T}(\Dcal, K')$. That is, there is some topological order $\sigma \in \Tcal(\Dcal, K')$ such that $W_k \le W_\kappa$ only if  $\comps^{-1}(\sigma, \kappa) \le \comps^{-1}(\sigma, k)$. 
\end{cor}
{\sc Proof}: See Appendix \ref{sec:proof_discussion}. $\Box$

\cref{cor:implementable} provides a necessary condition for waiting times to be implementable. While completely characterising the set of implementable waiting times for a particular $\Lambda$, $\mu$, and $M$ is difficult in general, we are able to provide a sufficient condition for waiting times to be implementable for menus such that the DAG satisfies the following property. 

\begin{dfn}\label{def:chain}{\rm (Chained  DAGs)} A DAG on $\mathbb{C}=\{\mathbb{C}_1, \mathbb{C}_2,\dots,\mathbb{C}_{K}\}$ is {\rm chained} if there exists a partition $\mathscr{C}=\{\mathscr{C}_1,\mathscr{C}_2,\dots,\mathscr{C}_L\}$ of $\mathbb{C}$ such that the DAG includes a directed arc from $\mathbb{C}_i$ to $\mathbb{C}_k$ if and only if $\mathbb{C}_i \in \mathscr{C}_\ell$ and $\mathbb{C}_k \in \mathscr{C}_{\ell+1}$ for some $\ell \in [L-1]$.
\end{dfn}

 \cref{fig:Ex_chained_DAG} illustrates an example of a chained DAG in panel (a) and one unchained DAG
 (i.e., a DAG that is not chained) in panel (b), both over a collection of seven CRP components. 
\begin{figure}[h] \begin{center}
    \includegraphics[width=10cm]{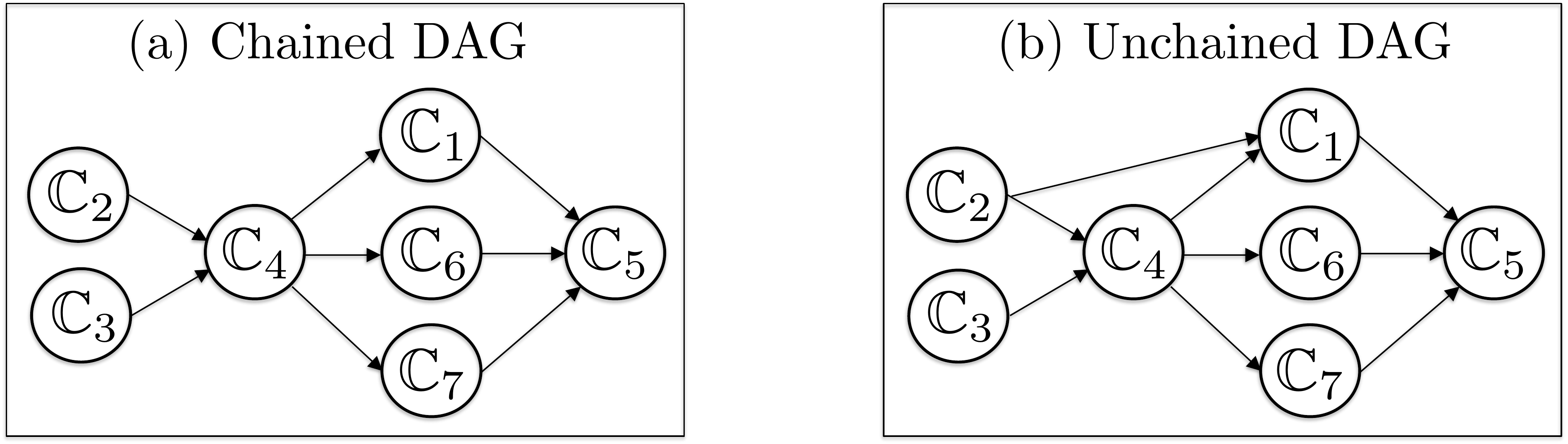}
    \caption[Model Description]{\sf \footnotesize Examples of chained (panel a) and unchained (panel b) DAGs over seven CRP components.}
     \label{fig:Ex_chained_DAG}
    \end{center}
\end{figure}
For the chained DAG in panel (a), $L=4$ and $\mathscr{C}_1=\{\mathbb{C}_2,\mathbb{C}_3\}$, $\mathscr{C}_2=\{\mathbb{C}_4\}$, $\mathscr{C}_3=\{\mathbb{C}_1, \mathbb{C}_6, \mathbb{C}_7\}$ and $\mathscr{C}_4=\{\mathbb{C}_5\}$. On the other hand, to see that the DAG in panel (b) is not chained, note that we cannot satisfy the requirement in \cref{def:chain} if we consider the three CRP components $\mathbb{C}_1$, $\mathbb{C}_2$ and $\mathbb{C}_4$. Indeed, the arcs connecting $\mathbb{C}_2$ and $\mathbb{C}_4$ to $\mathbb{C}_1$ would require that $\mathbb{C}_2$ and $\mathbb{C}_4$ belong to the same class  $\mathscr{C}_l$ in the partition $\mathscr{C}$ for some $\ell$, but then the arc connecting $\mathbb{C}_2$ to $\mathbb{C}_4$ would require these two CRP components to be in different classes in $\mathscr{C}$.

For menus such that the DAG is chained, the following result regarding which vectors of waiting times are implementable applies. 

\begin{prop}\label{prop:Wchained} Take limiting arrival rates $\Lambda$, service rates $\mu$, and a menu $M$ such that $M \in \mathcal{M}_+(\Lambda, \mu)$, and the collection of CRP components $\mathbb{C}=\{\mathbb{C}_1,\mathbb{C}_2,\dots,\mathbb{C}_K\}$ and the chained DAG $\Dcal = (|K|, \mathcal{A})$ are induced. Let $\mathscr{C}=\{\mathscr{C}_1,\mathscr{C}_2,\dots,\mathscr{C}_L\}$ be the partition of $\mathbb{C}$ described in \cref{def:chain}. 

The vector $W=(W_1,W_2,\dots,W_K) \in \RR^K_+$ is implementable if the following both hold:
\begin{enumerate}[nosep]
    \item[\rm (i)] $W_k = W_\kappa$ for all $(k, \kappa) \in [K] \times [K]$ such that $W_k \in \mathscr{C}_\ell$ and $W_\kappa \in \mathscr{C}_\ell$ for some $\ell \in [L]$,
    \item[\rm (ii)]  $W_k < W_\kappa$ for all $(k, \kappa) \in [K] \times [K]$ such that $W_k \in \mathscr{C}_\ell$ and $W_\kappa \in \mathscr{C}_{\ell'}$ for some $(\ell, \ell') \in [L] \times [L]$ where $\ell < \ell'$. 
\end{enumerate}
\end{prop}
{\sc Proof}: See Appendix \ref{sec:proof_discussion}. $\Box$

This tells us that we greatly increase the set of implementable outcomes by using a more general heavy traffic scaling.  

\subsection{Menu Design}
We now turn our attention to some simple questions regarding the design of menus of customer classes. We will consider two objectives: (1) minimising the total average delay across all customer classes, and (2) minimising the maximum expected delay of any customer class. We will assume that the arrival rates into the customer classes $\htp{\lambda}$ and the service rates $\mu$ are fixed, and the service provider is designing the menu $M$, or the compatibility between the customer classes and servers. 

When the service provider has complete flexibility over how to design the menu, the service provider can minimise both the average delay and the maximum delay faced by any customer class simultaneously. The following proposition shows that this can be achieved with a menu that has a single CRP component. 
\begin{prop}\label{prop:min_avg_waits}
Given arrival rates $\htp{\lambda}$ and service rates $\mu$ satisfying \cref{assum:heavy_traffic}, for any admissible menu $M \in \mathcal{M}(\htp{\lambda}, \mu)$,
\begin{equation*}
        \widehat{W}_{\mathbb{C}_k} \ge \frac{1}{|\Gamma|},
\end{equation*}
for all $k \in [K]$.

Furthermore, $\widehat{W}_{\mathbb{C}_k} = \frac{1}{|\Gamma|}$ for some $k \in [K]$ if and only if there exists a directed path from $\widehat{W}_{\mathbb{C}_k}$ to any other CRP component $\mathbb{C}_\kappa$ with $\kappa \in \{ [K]\setminus k\}$. This condition is trivially satisfied if there is only one CRP component. 
\end{prop}
{\sc Proof}: See Appendix \ref{sec:proof_discussion}. $\Box$

Therefore any menu that induces a single CRP component will ensure that all customer classes achieve the minimum possible expected delay, hence minimising both the average delay across all customer classes and the maximum delay faced by any customer class. The following proposition is helpful in designing such a menu. 
\begin{prop}\label{prop:singleCRP_condition}
Consider a system with limiting arrival rates $\Lambda$ and service rates $\mu$. Any menu $M$ such that 
\begin{equation*}
    \sum_{j \in \mathscr{S}}\sum_{i \in [n]} \Lambda_i m_{ij} < \sum_{j \in \mathscr{S}} \mu_j, \quad \text{for all } \mathscr{S} \subsetneq [m]
\end{equation*}
will be admissible for any vector of slacks $\Gamma \in \RR^n$ such that $|\Gamma| > 0$. Furthermore, such a menu will induce a single CRP component. 
\end{prop}
{\sc Proof}: See Appendix \ref{sec:proof_discussion}. $\Box$

A complete menu, in which every customer class is compatible with every server, will always satisfy this condition. The complete menu will operate like a single queue served by all servers according to an FCFS service discipline. \cref{prop:singleCRP_condition} also tells us that we do not need to know the values for the slacks $\Gamma$ to design a delay minimising menu, making it easier to implement in practice. 

While a menu that induces a single CRP component minimises delays, it may not be desirable or even feasible to offer such a menu due to real-world compatibility constraints on which servers can serve which customer types. Motivated by these sorts of constraints, we consider the question of how to design the DAG on a collection of CRP components to minimise expected delays for customers. 

It will be useful first to understand the expression for average expected delays across all customer classes. In \cref{eq:WCRP_def} we defined the delay of each CRP component conditional on being in a particular topological order. We can similarly define $\bar{w}_{\sigma}$, the average delay across all customer classes conditional on being in a particular topological order $\sigma$, as
\begin{equation}\label{eq:average_conditional_delay}
     \bar{w}_{\sigma}= \sum_{\kappa = 1}^{K'} \frac{\sum_{k = 1}^{\kappa} \widetilde{\mu}_{\sigma(k)}}{\sum_{\ell = 1}^{\kappa}\widetilde{\gamma}_{\comps(\sigma, \ell)}}. 
\end{equation}

This then lets us express the average expected delay for a particular menu $M$ as 
\begin{align}
\label{eq:average_wait}
\bar{W} = \frac{1}{|\mu|}\sum_{t=1}^{T(M)} \left({\mathbb{Q}(\sigma_t) \over \mathbb{Q}(\sigma_1)+\mathbb{Q}(\sigma_2)+\cdots+\mathbb{Q}(\sigma_{T(M)})}\right)\, \sum_{\kappa = 1}^{K'} \frac{\sum_{k = 1}^{\kappa} \widetilde{\mu}_{\sigma(k)}}{\sum_{\ell = 1}^{\kappa}\widetilde{\gamma}_{\comps(\sigma_t, \ell)}}. 
\end{align}
Here we can also see the differences with \cite{Afecheetal2019}, in which the authors find that the average delays depend only on the number of CRP components. With our more general scaling, the average delays depend on the values of the slacks themselves, as well as the structure of the DAG and the set of topological orders that are induced. 

Introducing additional arcs into the DAG reduces the number of topological orders. If we can introduce or remove arcs from a DAG in such a way that the system spends more time in states associated with topological orders that have lower conditional average delays $\bar{w}_{\sigma}$, then the total average delay will be reduced. However, the values of the slacks of the different CRP components $\widetilde{\gamma}$ limit how we are able to adjust the DAG and still have an admissible menu. This leads us to the following definition of an admissible topological order.  

\begin{dfn} A topological order $\sigma$ is \emph{admissible} for arrival rates $\htp{\lambda}$ and service rates $\mu$ satisfying \cref{assum:heavy_traffic}, and a collection of CRP components $\{ \mathbb{C}_1,\dots,\mathbb{C}_{K'}, \mathbb{C}_{K'+1}, \dots, \mathbb{C}_{K} \}$ if if $\sum_{\ell = 1}^{k} \widetilde{\gamma}_\ell > 0$ for all $k \in [K']$.  
\end{dfn}

The following lemma tells us how admissible topological orders relate to admissible menus.  
\begin{lem}\label{lem:admissible_top}
Take any arrival rates $\htp{\lambda}$ and service rates $\mu$ satisfying \cref{assum:heavy_traffic}, and any collection of CRP components $\{ \mathbb{C}_1,\dots,\mathbb{C}_{K'}, \mathbb{C}_{K'+1}, \dots, \mathbb{C}_{K} \}$. For any admissible topological order $\sigma$, we can construct an admissible menu $M \in \mathcal{M}(\htp{\lambda}, \mu)$ such that the DAG induced by $M$ with $\htp{\lambda}$ and $\mu$ only admits the topological order $\sigma$. Furthermore, if $\sigma$ is not admissible, then there are no admissible menus $M$ that admit the topological order $\sigma$. 
\end{lem}
{\sc Proof}: See Appendix \ref{sec:proof_discussion}. $\Box$

The set of admissible topological orders tells us which DAGs are feasible given a particular CRP component. We can then minimise average delays by identifying the topological order with the lowest condition delays. 
\begin{prop}\label{prop:delay_min_DAG}
Given limiting arrival rates $\Lambda$, service rates $\mu$, slacks $\Gamma$, and CRP components $\{\mathbb{C}_1,\dots,\mathbb{C}_{K'}, \mathbb{C}_{K'+1}, \dots, \mathbb{C}_{K} \}$, there will be a permutation of CRP components $\sigma$ that minimises the average expected delay across all implementable topological orders,
\begin{equation*}
     \bar{w}_{\sigma}= \sum_{\kappa = 1}^{K'} \frac{\sum_{k = 1}^{\kappa} \widetilde{\mu}_{\sigma(k)}}{\sum_{\ell = 1}^{\kappa}\widetilde{\gamma}_{\comps(\sigma, \ell)}}.
\end{equation*}
The DAG or menu that will minimise delays is one that only allows for this topological order. 
\end{prop}
{\sc Proof}: See Appendix \ref{sec:proof_discussion}. $\Box$

Given that adding arcs to a DAG is achieved by adding additional flexibility to a service system, one might think that adding an additional arc to a DAG will always reduce expected delays. However, we find that adding arcs to the DAG may potentially increase, decrease, or not affect the average delays. This can be shown through the following two server example. 

Consider the case of two independent $M/M/1$ queues. We will use $M_a$ to denote this menu. Let the arrivals rates be $\htp{\lambda}_1 = 1 - \eps\gamma_1$, and $\htp{\lambda}_2 = 1 - \eps\gamma_2$, and let $\mu_1 = \mu_2 = 1$. It is straightforward to calculate that $\widehat{W}_1 = 1/\gamma_1$ and $\widehat{W}_2 = 1/\gamma_2$. The average delay across both customer classes is then
\begin{equation}
    \bar{W}_a = \frac{1}{2} \left(\frac{1}{\gamma_1} + \frac{1}{\gamma_2} \right)
\end{equation}

If we were to consider the alternative menu 
\begin{equation}
M_b=\left[\begin{array}{cc} 1 & 1 \\ 0 & 1\end{array}\right],
\end{equation}
then using \cref{thm:CRPdelay_append} we find that $\widehat{W}_1 = 1/(\gamma_1 + \gamma_2)$ and $\widehat{W}_2 = 1/(\gamma_1 + \gamma_2) + 1/\gamma_2$. The average delay across both customer classes is then
\begin{equation}
    \bar{W}_b = \frac{1}{\gamma_1 + \gamma_2} + \frac{1}{2\gamma_2}. 
\end{equation}
Therefore the difference in average delays is 
\begin{equation*}
    \Delta_{ab} :=  \bar{W}_b-  \bar{W}_a= \frac{1}{\gamma_1 + \gamma_2} - \frac{1}{2\gamma_1}.
\end{equation*}
When $\gamma_1 = \gamma_2$, $\Delta_{ab} = 0$ and menus $M_a$ and $M_b$ have the same average delays. When $\gamma_1 > \gamma_2$, $\Delta_{ab}$ is positive, and menu $M_b$ has higher average delays than $M_a$, despite the additional flexibility. Otherwise, $\Delta_{ab}$ is negative, and menu $M_b$ has lower average delays than $M_a$. 

This simple example demonstrates that adding additional flexibility to the design of the menu does not necessarily reduce the average delay (i.e., some form of Braess's paradox). Therefore if a service provider is considering adding additional flexibility to a system, it is important to carefully consider the way in which flexibility is being added. 

\subsection{Numerical example}\label{sec:numex}
We will end this section by returning to our example in \cref{fig:residual_matching} (a) to make some of the ideas discussed in the section more concrete. Recall the menu $M$ is given by
\begin{equation}\label{eq:discussion_menu}
M=\left[\begin{array}{cccc} 1 & 0 & 0 & 0 \\ 0 & 1 & 0 & 0 \\ 0 & 0 & 1 & 0 \\ 1 & 1 & 1 & 1\end{array}\right].
\end{equation}
The limiting arrival rates are $\Lambda = (2, 1, 1, 2)$, and service rates are $\mu = (2, 1, 2, 1)$. We will let the sequence of arrival rates be $\htp{\lambda}_i = \Lambda_i - \epsilon\gamma_i$ for $1 \le i \le 4$. We have three CRP components, $\mathbb{C}_1$ consisting of class 1 and server 1, $\mathbb{C}_2$ consisting of class 2 and server 2, and $\mathbb{C}_3$ consisting of classes 3 and 4 and servers 3 and 4. 
 
We will begin by considering the question of implementability. We can see that the DAG induced by $M$ is a chained DAG, with $\mathbb{C}_1$ and $\mathbb{C}_2$ belonging to one partition in the chain, and $\mathbb{C}_3$ belonging to the other partition in the chain. Then \cref{prop:Wchained} tells us that we can implement any waiting times $W_1 = W_2 > W_3 > 0$. 

In this simple case, we can see which delays are implementable more directly, by looking at the exact expressions for the delays. Using \cref{thm:CRPdelay_append}, we can calculate the delays as 
\begin{equation*}
    \hat{W}_1 = \frac{1}{{\gamma}_1} + \frac{1}{{\gamma}_1 + {\gamma}_2 + {\gamma}_3 + \gamma_4}, \quad \hat{W}_2 = \frac{1}{{\gamma}_2} + \frac{1}{{\gamma}_1 + {\gamma}_2 + {\gamma}_3}, \quad \text{and } \hat{W}_3 = \frac{1}{{\gamma}_1 + {\gamma}_2 + {\gamma}_3 + \gamma_4}.
\end{equation*}
By looking at these expressions, we can see that we can implement any delays $W_1$, $W_2$, and $W_3$ such that $W_3 > 0$, $W_1 > W_3$ and $W_2 > W_3$. To do this we would let ${\gamma_1} = \frac{1}{W_1 - W_3}$, ${\gamma_2} = \frac{1}{W_2 - W_3}$, and $\gamma_3  + \gamma_4 = \gamma_1 + \gamma_2 - 1/W_1$. 

This also suggests that in a congested system, a service provider is able to produce significant improvements in delay if they can make small changes to the arrival rates into the different service classes. 

Suppose arrival rates are initially such that the slacks are proportional to arrival rates, i.e. $\gamma = \Lambda$, as in \cite{Afecheetal2019}. The following table shows us the improvements in delay by adjusting the slacks so that $\gamma' = (9, 9, -3, -9)$ for different values of $\epsilon$. Note that $|\Lambda| = |\gamma'|$, so this adjustment does not alter the total arrival rate of customers into the system. We also show the percentage difference in average delays, denoted $\delta\bar{W}\%$, as well as the percentage of customers who are joining a different customer class across the two scenarios, denoted $\delta\lambda\%$. 

\begin{table}[h]
\centering
\begin{tabular}{||c|c|c|c|c|c|c|c|c||}
\hline \hline
$\epsilon$            & $\gamma$  & $\hat{W}_1$ & $\hat{W}_2$ & $\hat{W}_3$ & $\hat{W}_4$ & $\bar{W}$ & $\delta\bar{W}\%$               & $\delta\lambda\%$             \\ \hline \hline
\multirow{2}{*}{0.1}  & $\Lambda$ & 0.5727      & 1.0652      & 0.1649      & 0.1182      & 0.4353    & \multirow{2}{*}{61.43\%} & \multirow{2}{*}{33.33\%} \\ \cline{2-7}
                      & $\gamma'$ & 0.2029      & 0.2001      & 0.2344      & 0.1237      & 0.1679    &                         &                         \\ \hline \hline
\multirow{2}{*}{0.05} & $\Lambda$ & 0.6171      & 1.1151      & 0.1651      & 0.1408      & 0.4660    & \multirow{2}{*}{60.32\%} & \multirow{2}{*}{15.79\%} \\ \cline{2-7}
                      & $\gamma'$ & 0.2351      & 0.2339      & 0.1830      & 0.1431      & 0.1849    &                         &                         \\ \hline \hline
\multirow{2}{*}{0.01} & $\Lambda$ & 0.6563      & 1.1562      & 0.1662      & 0.1612      & 0.4929    & \multirow{2}{*}{56.82\%} & \multirow{2}{*}{3.03\%} \\ \cline{2-7}
                      & $\gamma'$ & 0.2678      & 0.2677      & 0.1670      & 0.1613      & 0.2128    &                         &                         \\ \hline \hline
\end{tabular}
\end{table}

As we can see, significant improvements in scaled delays are achieved while only changing the arrivals of a relatively small fraction of customers, with the improvements in comparison to the change required increasing as congestion increases. 

Finally, we look at the question of menu design. In particular, we look at how we can change a menu to improve delays given a fixed CRP component structure, and fixed arrival rates. The residual matching for the menu $M$ in \cref{eq:discussion_menu} with limiting arrival rates $\Lambda = (2, 1, 2, 1)$ and service rates $\mu = (2, 1, 1, 2)$ is 
\begin{equation}\label{eq:discussion_menu}
\breve{M}=\left[\begin{array}{cccc} 1 & 0 & 0 & 0 \\ 0 & 1 & 0 & 0 \\ 0 & 0 & 1 & 0 \\ 0 & 0 & 1 & 1\end{array}\right].
\end{equation}
There are 6 possible permutations of CRP components when the menu is just the residual matching $\breve{M}$, these permutations being all the permutations of the number (1, 2, 3). We can use \cref{eq:average_conditional_delay} to calculate the expected delay conditional on a particular permutation of CRP components. In this case, we will assume the slacks are $\gamma = (4, 3, 1, 1)$. The following table uses \cref{eq:average_conditional_delay} to calculate the conditional delays for all possible server permutations. 
\begin{table}[h]
\centering
\begin{tabular}{|c|c|}
\hline
Permutation & Delay \\ \hline
(1,2,3)   & $\frac{\mu_1}{{\gamma}_1} + \frac{\mu_1 + \mu_2}{\gamma_1 + \gamma_2} + \frac{\mu_1 + \mu_2 + \mu_3 + \mu_4}{\gamma_1 + \gamma_2 + \gamma_3 + \gamma_4} = 1.595$     \\ \hline
(1,3,2)     & $\frac{\mu_1}{{\gamma}_1} + \frac{\mu_1 + \mu_3 + \mu_4}{\gamma_1 + \gamma_3 + \gamma_4} + \frac{\mu_1 + \mu_2 + \mu_3 + \mu_4}{\gamma_1 + \gamma_2 + \gamma_3 + \gamma_4} = 2$      \\ \hline
(2,1,3)     & $\frac{\mu_2}{{\gamma}_2} + \frac{\mu_1 + \mu_2}{\gamma_1 + \gamma_2} + \frac{\mu_1 + \mu_2 + \mu_3 + \mu_4}{\gamma_1 + \gamma_2 + \gamma_3 + \gamma_4} = 1.429$      \\ \hline
(2,3,1)     &  $\frac{\mu_2}{{\gamma}_2} + \frac{\mu_2 + \mu_3 + \mu_4}{\gamma_2 + \gamma_3 + \gamma_4} + \frac{\mu_1 + \mu_2 + \mu_3 + \mu_4}{\gamma_1 + \gamma_2 + \gamma_3 + \gamma_4} = 8$     \\ \hline
(3,1,2)     & $\frac{\mu_3 + \mu_4}{{\gamma}_3 + \gamma_4} + \frac{\mu_1 + \mu_3 + \mu_4}{\gamma_1 + \gamma_3 + \gamma_4} + \frac{\mu_1 + \mu_2 + \mu_3 + \mu_4}{\gamma_1 + \gamma_2 + \gamma_3 + \gamma_4} = 3$      \\ \hline
(3,2,1)     & $\frac{\mu_3 + \mu_4}{{\gamma}_3 + \gamma_4} + \frac{\mu_2 + \mu_3 + \mu_4}{\gamma_2 + \gamma_3 + \gamma_4} + \frac{\mu_1 + \mu_2 + \mu_3 + \mu_4}{\gamma_1 + \gamma_2 + \gamma_3 + \gamma_4} = 2.967$   \\ \hline
\end{tabular}
\end{table}
We can see from this table that the permutation of CRP components that minimises delay is (2,1,3). We can then design a menu such that the DAG only admits this specific topological order. The DAG that achieves this is shown below. 

\begin{figure}[h] \begin{center}
    \includegraphics[width=10cm]{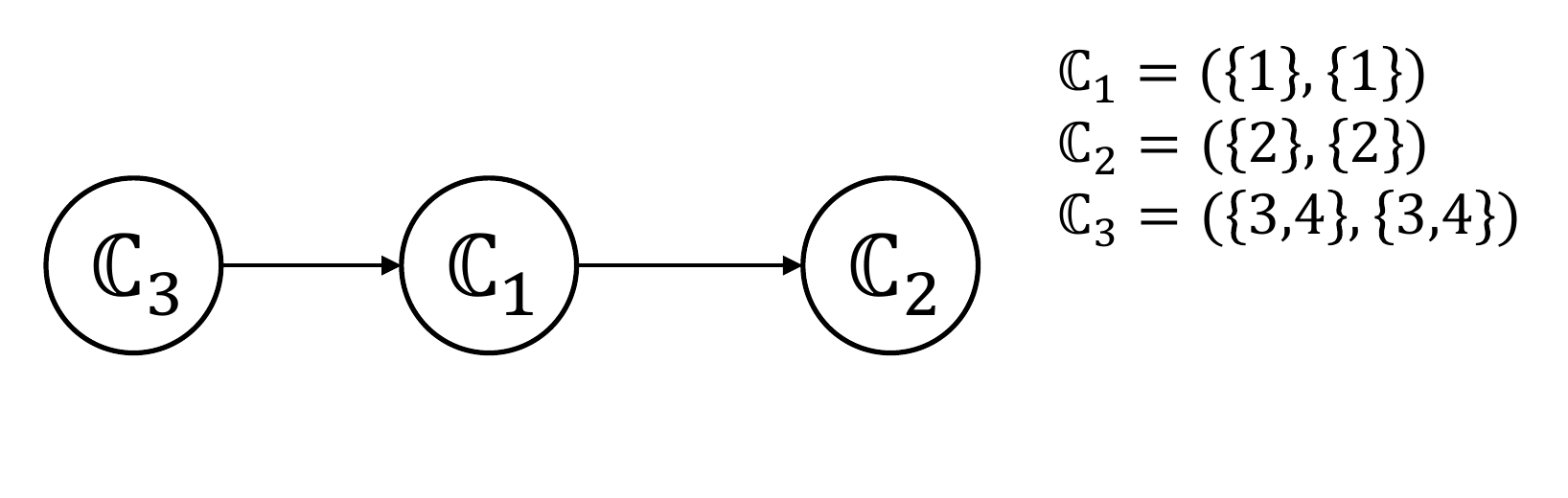}
    \caption[Model Description]{\sf \footnotesize Delay minimising DAG.}
     \label{fig:Ex_chained_DAG}
    \end{center}
\end{figure}

This DAG can by having the customer class in $\mathbb{C}_1$ served by the server in $\mathbb{C}_2$, and either of the customer classes in $\mathbb{C}_3$ served by the server in $\mathbb{C}_1$. The following menu is one example of a menu that achieves this.   
\begin{equation}\label{eq:discussion_menu}
{M'}=\left[\begin{array}{cccc} 1 & 0 & 0 & 0 \\ 1 & 1 & 0 & 0 \\ 0 & 1 & 1 & 0 \\ 0 & 0 & 1 & 1\end{array}\right].
\end{equation}

In comparison, the original menu $M$ in \cref{eq:discussion_menu} with $\gamma = (4, 2, 1, 1)$ has average delays of 1.5, which as expected is higher than the average delays of our newly designed menu.

\section{Proof of Main Results}\label{sec:ProofTheorem}

\subsection{Proof of \cref{thm:CRPdelay_append}}\label{sec:waits_proof}
The key observation needed to prove \cref{thm:CRPdelay_append} is that only a relatively small subset of states have positive probability in heavy-traffic, and the information about which states have positive probability is captured by the CRP components and the DAG on the CRP components. However, before we go into more detail, it will be useful to introduce some notation. In section \cref{eq:CRP_mu_lambda}, we defined the aggregate arrival rate for a CRP component $\mathbb{C}_k$ to be  $\htp{\widetilde{\lambda}_k} = \sum_{i \in \Ccal_k} \htp{\lambda_i} = \widetilde{\Lambda}_k - \epsilon \widetilde{\gamma}_k + o(\epsilon)$. For a subset of servers $S \subseteq [m]$, we define the {\it slack} for $S$ by:
\begin{align}\label{def:capacity_slack}
    \Delta(S) &= \mu_S - \lambda_{U_S(M)},
\end{align}
where $U_S(M)$ is defined in Proposition~\ref{prop:stability} as the subset of service classes that can only be served (or, uniquely served) by servers in $S$ under the menu $M$. For succinctness, we will suppress the dependence on $M$ in this section and use the notation $U(S)$ for $U_S(M)$.

It will also be useful to further aggregate the state space described in \cref{sec:adanweissresult} so that the state depends only on the server permutation $s$ and the number of busy servers $b$, and not the number of customers. Specifically, for a server permutation $s = \{s_1,\ldots, s_m\}$ and $b \in \{0,1,\ldots,m\}$ define:
\[ P(s;b) = \{ x \in X : x = (s_1, n_1,\ldots, s_b, n_b, s_{b+1},s_{b+2}\ldots, s_m)\} \]
as the set of all states where $s$ is the ranking of servers in terms of the age of the customer for busy servers and the time since idleness for idle servers, and where exactly the first $b$ servers in $s$ are busy. We then have the following expression for the probability of the aggregate state $P(s;b)$:
\begin{align}
\pi(P(s;b)) &= \sum_{n_1 = 0}^{\infty} \cdots \sum_{n_b = 0}^ \infty {\cal B}  \prod_{\ell=1}^b \frac{\lambda^{n_\ell}_{U(s_1,\ldots, s_\ell)} }{\mu_{\{s_1,\ldots, s_\ell\}}^{n_\ell+1} } \prod_{\ell=b+1}^m \lambda_{C(s_\ell, \ldots, s_m)}^{-1} \nonumber \\
& =  {\cal B}\, \;\prod_{\ell=1}^b \frac{1} {\Delta(s_1,\ldots, s_\ell)} \prod_{\ell=b+1}^m \lambda_{C(s_\ell, \ldots, s_m)}^{-1}. \label{eq:steadyP_k}
\end{align}

As a last step before developing the proof of \cref{thm:CRPdelay_append}, in \cref{lem:CRP_properties} we state some properties of CRP components and topological orders that will be useful. This lemma has been slightly modified from \cite[Lemma 6]{Afecheetal2019}.  
\begin{lem} \label{lem:CRP_properties} Let $M$ be a service menu and $\{ \mathbb{C}_1, \ldots, \mathbb{C}_{K'}, \mathbb{C}_{K'+1}, \ldots, \mathbb{C}_K \}$ be its CRP components under a given heavy-traffic equilibrium strategy profile. For a CRP component $\mathbb{C}_k = (\Ccal_k, \Scal_k)$ with non-empty $\Scal_k$ (i.e., $k \in [K']$) :
\begin{enumerate}
\item[\rm (i)] The aggregate demand of service classes converges to the aggregate service rate as $\epsilon \to 0$, that is, $\widetilde{\Lambda}_k := \Lambda_{\Ccal_k} = \mu_{\Scal_k} =: \widetilde{\mu}_k$ (see \eqref{eq:CRP_mu_lambda} for definitions).

\item[\rm (ii)] For any strict subset of servers $S \subset \Scal_k$, the set of service classes in residual matching $\breve{M}$ served only by $S$ is a strict subset of $\Ccal_k$, and $S$ exhibits strictly positive slack as $\epsilon \to 0$, that is, 
\[ \forall S \subset \Scal_k: U_S(\breve{M}) \subset \Ccal_k  \quad \mbox{and} \quad  \mu_{S} > {\Lambda}_{U_S(\breve{M})}.\]
Further, since $U_S(M) \subseteq U_S(\breve{M})$, the positive slack condition also holds for $U_S(M)$.

(Recall that $U_S(M)$ is the subset of service classes that can only be served by servers in $S$.)
\end{enumerate}

Let $\sigma \in \Tcal(\Dcal, K')$ be a topological order of the CRP components with non-empty server sets. Define $\mathscr{S}_k = \Scal_{\sigma(1)} \cup \Scal_{\sigma(2)} \cup \cdots \cup \Scal_{\sigma(k)}$ and $\mathscr{C}_k= \Ccal_{\sigma(1)} \cup \Ccal_{\sigma(2)} \cup \cdots \cup \Ccal_{\sigma(k)}$ to be the subset of servers and customer classes in the first $k$ CRP components in the topological order. Define 
\[ \mathscr{C}'_k = \left\{ \cup_{\kappa} \Ccal_{\kappa} |  \kappa \in \{K'+1,\ldots, K\} : \exists k' \in \{ 1, \ldots, k\}, \kappa  \in \comps(\sigma, k') \right\}\]
to be the service classes of server-less CRP components that are part of $\comps(\sigma,k')$ for some $k' \in [k]$. Then,

\begin{enumerate}

\item[\rm (iii)] Customers  in $\mathscr{C}_k \cup \mathscr{C}'_k$ are exclusively served by servers in $\mathscr{S}_k$. That is,
$$U_{\mathscr{S}_k}(M) = \mathscr{C}_k \cup \mathscr{C}'_k .$$

\item[\rm (iv)] The capacity slack of the set of servers $\mathscr{S}_k$ converges to zero as $\epsilon \to 0$, in particular, 
$$\Delta(\mathscr{S}_k) = \epsilon \sum_{\ell = 1}^k \widetilde{\gamma}_{\comps(\sigma, \ell)} + o(\epsilon).$$
\end{enumerate}
\end{lem}
{\sc Proof}: See Appendix \ref{sec:proof_waits_appen}. $\Box$

We can now begin calculating the expected waits. Using the aggregated states from \cref{eq:steadyP_k}, the following lemma (rephrased) from \cite{Afecheetal2019} gives an expression for the mean waiting time for each service class in terms of the probabilities $\pi(P(s;b))$.
\begin{lem}{\cite[Lemma 6]{Afecheetal2019}} \label{lem:WaitingTime} The steady-state mean waiting time of service class $i$ is equal to
 $$W_i = \sum_{s \in \Sigma_m} \sum_{b=1}^m W_i(s;b) \cdot \pi(P(s;b)),$$
 where $\Sigma_m$ denotes the set of all the permutations of $[m]$, 
 \[ W_i(s;b) = \sum_{\ell=1}^b {\ind\big(i \in U(s_1,\ldots, s_\ell)\big) \over \Delta(s_1,\dots,s_\ell)}, \]
 and $\pi(P(s;b))$ is given by \eqref{eq:steadyP_k}.
\end{lem}

We are able to simplify these expressions further by showing that only a relatively small subset of aggregate states $(s, b)$ have asymptotically non-zero probabilities in heavy-traffic. These states are exactly those that are consistent with $\Tcal(\Dcal, K') = (\sigma_1, \ldots, \sigma_T)$ the collection of topological orders on $\{\mathbb{C}_1, \ldots, \mathbb{C}_{K'}\}$, a notion we will formalize in \cref{def:permutationinduced}. Our first step to showing this is to consider the slacks $\Delta(s_1,\ldots, s_\ell)$, which the preceding lemma suggests will be an important part of the analysis. \cref{lem:criticalservers} below, which is an extension of \cite[Lemma 4]{Afecheetal2019} shows that only certain subsets of servers have ``interesting'' slacks under a given sequence of arrival rates $\htp{\lambda}$.

\begin{lem}\label{lem:criticalservers} Let $\Dcal$ be the DAG for the CRP decomposition $\{\mathbb{C}_1,\ldots, \mathbb{C}_{K'},\mathbb{C}_{K'+1}, \ldots, \mathbb{C}_K\}$ under some menu $M$ and a given heavy-traffic equilibrium strategy profile. Then, a subset of servers $\{s_1, \ldots, s_\ell\} \subseteq [m]$ satisfies 
$$\lim_{\epsilon \to 0} {\epsilon \over \Delta(s_1, \ldots, s_\ell)} > 0$$
if and only if there exists a topological order $\sigma \in \Tcal(\Dcal, K')$ and an integer $k$ such that 
\begin{equation}\label{eq:criticalservers}\{s_1,\dots,s_\ell\} = \bigcup_{i=1}^k \Scal_{\sigma(i)}.\end{equation}
Further, in this case :
\[ \lim_{\epsilon \to 0} {\epsilon \over \Delta(s_1, \ldots, s_\ell)} = \frac{1}{ \sum_{i = 1}^k \widetilde{\gamma}_{\comps(\sigma, i)}} \]
for any topological order $\sigma$ for which \eqref{eq:criticalservers} is satisfied.
\end{lem}
{\sc Proof}: See Appendix \ref{sec:proof_waits_appen}. $\Box$

As implied in the previous paragraph, we can use \cref{lem:criticalservers} to prove \cref{prop:steadystate} below, which states that a relatively small number of aggregate states have positive steady-state probability in heavy traffic; these are the aggregate states $P(s;m)$ in which $s$ is a permutation of the servers induced by a topological order $\sigma \in \Tcal(\Dcal, K')$ and such that all servers are busy.

\begin{dfn}\label{def:permutationinduced} {\em ({\sf Server Permutations Induced by Topological Orders})}   We say that a permutation of the servers $s=(s_1,s_2,\dots,s_m) \in \Sigma_m$ is {\em induced by} the topological order $\sigma \in \Tcal(\Dcal, K')$, if $s$ can be expressed as a concatenation of sub-permutations:
\[ s = \left( \sbb_{\sigma(1)} || \sbb_{\sigma(2)}|| \cdots || \sbb_{\sigma(K')} \right) \]
with $\sbb_k \in \Sigma_{\Scal_k}$ denoting a permutation of the servers $\Scal_{k}$ of CRP component $\mathbb{C}_k$. In other words, the servers of a CRP component are contiguous in the permutation $s$, and the order of the CRP components obeys the topological order $\sigma$.

\end{dfn}

Returning to our four server example in \cref{subfig:DAG_a}, the CRP components were $\mathbb{C}_1 = (\mathcal{C}_1, \mathcal{S}_2 = (\{1\}, \{1\})$, $\mathbb{C}_2 = (\mathcal{C}_2, \mathcal{S}_2 = (\{2\}, \{2\})$, and $\mathbb{C}_3 = (\mathcal{C}_3, \mathcal{S}_3 = (\{3, 4\}, \{3, 4\})$, and the topological orders were $\sigma_a = (1, 2, 3)$ and $\sigma_b = (2, 1, 3)$. \cref{def:permutationinduced} tells us the topological order $\sigma_a$ induces two possible server permutations, $s_{a1} = (s_1||s_2||s_3||s_4)$ and $s_{a2} = (s_1||s_2||s_4||s_3)$. 

The next proposition is an extension of \cite[Proposition 2]{Afecheetal2019}.

\begin{prop}\label{prop:steadystate} 
 Let $\Dcal$ be the DAG for the CRP decomposition $\{\mathbb{C}_1,\ldots, \mathbb{C}_{K'},\mathbb{C}_{K'+1}, \ldots, \mathbb{C}_K\}$ under some menu $M$ and a heavy-traffic strategy profile. Let $s \in \Sigma_m$ be a server permutation.   
 \begin{enumerate}
     \item If $b < m$, and/or $s$ is not a permutation of the servers induced by some topological order $\sigma \in \Tcal(\Dcal, K')$, then
     \[ \lim_{\epsilon \to 0} \pi(P(s;b)) =0. \]
     \item If $b=m$ and $s = \left( \sbb_{\sigma(1)} || \sbb_{\sigma(2)}|| \cdots || \sbb_{\sigma(K')} \right) $ is a server permutation induced by topological order $\sigma \in \Tcal(\Dcal, K')$ with subpermutations $\sbb_k \in \Sigma_{\Scal_k}$, then
          \[ \lim_{\epsilon \to 0} \pi(P(s;b)) = \Bcal' \cdot \mathbb{Q}(\sigma) \prod_{k=1}^{K'} \theta_k( \sbb_k) \]
          where $\Bcal'$ is a normalization constant, $\mathbb{Q}(\sigma)$ was defined in \eqref{eq:Q_def} as 
\[        \mathbb{Q}(\sigma) = \prod_{\kappa \in [K']} \frac{1}{\sum_{\ell=1}^\kappa  \widetilde{\gamma}_{\comps(\sigma, \ell)}},
\]
and $\left\{ \theta_k : \Sigma_{\Scal_k} \to \Re^+ \right\}_{k \in [K']}$ is a fixed collection of functions mapping the sub-permutation of servers of CRP components to positive reals.
 \end{enumerate}
\end{prop}

Using Proposition~\ref{prop:steadystate} and the normalization condition $\sum_{s \in \Sigma_m, 0 \leq  b \leq m } \pi(P(s;b))=1$, we get:
\begin{align*}
\lim_{\epsilon \to 0} \sum_{s \in \Sigma_m, 0 \leq  b \leq m } \pi(P(s;b)) & = \sum_{\sigma \in \Tcal(\Dcal,K')} \sum_{\substack{  s = (\sbb_{\sigma(1)} || \sbb_{\sigma(2)} || \cdots || \sbb_{\sigma(K')})  \\ \{\sbb_k \in \Sigma_{\Scal_k}\}_{k \in [K']}  } }  \pi(P(s;m)) \\ 
&= \sum_{\sigma \in \Tcal(\Dcal,K')} \sum_{\substack{  s = (\sbb_{\sigma(1)} || \sbb_{\sigma(2)} || \cdots || \sbb_{\sigma(K')})  \\ \{\sbb_k \in \Sigma_{\Scal_k}\}_{k \in [K']} } } \Bcal' \cdot \mathbb{Q}(\sigma) \prod_{k=1}^{K'} \theta_k(\sbb_k) \\
& = \left( \sum_{\sigma \in \Tcal(\Dcal,K')} \mathbb{Q}(\sigma) \right) \left( \Bcal' \sum_{ \{\sbb_k \in \Sigma_{\Scal_k}\}_{k \in [K']} } \prod_{k=1}^{K'} \theta_k(\sbb_k) \right) ,
\end{align*}
or,
\[ \left( \Bcal' \sum_{ \{\sbb_k \in \Sigma_{\Scal_k}\}_{k \in [K']} }   \prod_{k=1}^{K'} \theta_k(\sbb_k) \right) = \frac{1}{\sum_{\sigma \in \Tcal(\Dcal, K')} \mathbb{Q}(\sigma)}.\]

Finally, we provide a lemma giving expressions for the scaled $W_{i}(s;b)$ when $s$ is a server permutation induced by a topological order $\sigma$, and $b=m$, as these are the only permutations that will be important in arriving at the result. A somewhat remarkable fact is that the limiting scaled $W_{i}(s;m)$ depends only on the topological order $\sigma$ and not the full server permutation $s$.

\begin{lem}
\label{lem:scaled_W_serv_perm}
Let $s = (s_1, \ldots, s_m)$ be a server permutation induced by the topological order $\sigma \in \Tcal(\Dcal, [K'])$. For a service class $i \in \mathbb{C}_k$,
\begin{align}
    \lim_{\epsilon \to 0} \epsilon W_i(s;m) &= w_{\sigma, k} \ := \ \sum_{\kappa = \comps^{-1}(\sigma, k)}^{K'} \frac{1}{\sum_{\ell=1}^\kappa  \widetilde{\gamma}_{\comps(\sigma, \ell)} }.
\end{align}
\end{lem}
{\sc Proof}: See Appendix \ref{sec:proof_waits_appen}. $\Box$

Combining Proposition~\ref{prop:steadystate} with Lemmas~\ref{lem:WaitingTime}-\ref{lem:scaled_W_serv_perm}, the limiting scaled mean waiting time for service class $i \in \mathbb{C}_k$ is:
\begin{align*} 
    \widehat{W}_i^* &= \lim_{\epsilon \to 0} \epsilon \cdot W_i \\
    &= \lim_{\epsilon \to 0} \sum_{s \in \Sigma_m} \epsilon \sum_{b=1}^m W_i(s;b) \cdot \pi(P(s;b)). \\
    \intertext{Using the product rule of limits \footnote{Product rule of limits: If $\lim_{x \to x_0} f(x) = F$ and $\lim_{x\to x_0} g(x) = G$, then $\lim_{x \to x_0} f(x) g(x)$ exists and equals $FG$.} we can reduce the above sum to a sum over server permutations induced by topological orders, and where all servers are busy.}
    \widehat{W}_i^* &= \lim_{\epsilon \to 0} \sum_{\sigma \in \Tcal(\Dcal, K')} \ \ \sum_{\substack{  s = (\sbb_{\sigma(1)} || \sbb_{\sigma(2)} || \cdots || \sbb_{\sigma(K')})  \\ \{\sbb_k \in \Sigma_{\Scal_k}\}_{k \in [K']} } } \epsilon \cdot W_i(s;m) \cdot \pi(P(s;m)) \\ 
    &= \sum_{\sigma \in \Tcal(\Dcal, K')} \ \ \sum_{\substack{  s = (\sbb_{\sigma(1)} || \sbb_{\sigma(2)} || \cdots || \sbb_{\sigma(K')})  \\ \{\sbb_k \in \Sigma_{\Scal_k}\}_{k \in [K']} } }  w_{\sigma, k} \cdot \Bcal' \cdot \mathbb{Q}(\sigma) \prod_{\ell=1}^{K'}  \theta_\ell(\sbb_\ell) \\
    &= \sum_{\sigma \in \Tcal(\Dcal, K')}  w_{\sigma, k} \cdot \mathbb{Q}(\sigma) \ \ \sum_{\substack{  s = (\sbb_{\sigma(1)} || \sbb_{\sigma(2)} || \cdots || \sbb_{\sigma(K')})  \\ \{\sbb_k \in \Sigma_{\Scal_k}\}_{k \in [K']} } }   \Bcal'  \prod_{\ell=1}^{K'}  \theta_\ell(\sbb_\ell) \\
    &= \frac{\sum_{\sigma \in \Tcal(\Dcal, [K'])}  w_{\sigma, k} \cdot \mathbb{Q}(\sigma)}{\sum_{\sigma \in \Tcal(\Dcal, K')}   \mathbb{Q}(\sigma)}\\
    &=: \widetilde{W}_k,
\end{align*}
as in the theorem statement.

\subsection{Proof of \cref{thm:matching_independent_gamma}}\label{sec:proof_matching}
Throughout this section, we will take the menu $M$, limiting arrival rates $\Lambda$ and service rates $\mu$, and slacks $\Gamma$ to be given, and largely suppress any dependence on $M$ in the notation. We will let $\breve{M}$ be the residual matching of the menu $M$ with arrival rates $\Lambda$ and service rates $\mu$. 

Instead of directly working with the matching rates $\htp{p}_{ij}(M)$, we will look at the service probabilities $\htp{q}_{ij}$.  For all $i \in [n]$ and $j \in [m]$, $\htp{q}_{ij}(x)$ is the probability with which server $j$ serves customer $i$ given the system is in state $x$ and server $j$ has become idle. We prove \cref{thm:matching_independent_gamma} by deriving and simplifying expressions for the limiting service probabilities $q_{ij}$ for the menu $M$, and find that the limiting service probabilities depend only on the service rates $\mu$, limiting arrival rates $\Lambda$, and the connectivity within each CRP component. To do this, we will make use of a new state space aggregation which we will introduce here.

In \cref{sec:waits_proof}, we introduced the aggregate states $P(s,b)$ for ever $s \in \Sigma_m$ and $b \in [m]$. Recall that $P(s,b)$ is the set of all states where $s$ is the ranking of servers in terms of the age of the customers they are serving for busy servers, and the time since becoming idle for the idle servers, and $b$ is the number of busy servers. In this section, we further aggregate the state space, so that we can consider all of the states in which we observe a particular subpermutation of servers within a CRP component together. Specifically, for some $k \in [K']$ and some subpermutation $\sbb_k \in \Sigma_{\Sbb_k}$, we define 
\begin{equation*}
    P_k(\sbb_k) = \cup_{\sigma \in \Tcal(\Dcal, K')} \left\{ s \in P(s,m) | s = \left( \sbb_{\sigma(1)} || \cdots || \sbb_k || \cdots || \sbb_{\sigma(K')} \right), \sbb_{\kappa} \in \Sigma_{\Sbb_\kappa} \text{ for } \kappa \in [K'] \text{ and } \kappa \neq k\right\}
\end{equation*}
Note that while the set of aggregated states $P(s,b)$ does not depend on the menu being offered, $P_k(\sbb_k)$ depends on the set of topological orders, and hence does depend on the menu. 

The first main step of our derivation will be to calculate the limiting service probabilities for our new further aggregated state space. That is, for each pair of customer classes $i \in [n]$ and servers $j \in [m]$ in the same CRP component, and for any subpermutation of servers within that CRP component $\sbb_{k(k)} \in \Sigma_{\Sbb_{k(j)}}$, we would like to calculate $q_{ij}(P_{k(j)}(\sbb_{k(j)}))$, the limiting service probability of customer class $i$ by server $j$ given the system is in a state in $P_{k(j)}(\sbb_{k(j)})$. Recall that $k(j)$ denotes the index of the CRP component that server $j$ belongs to. We do not consider $i$ and $j$ that are not in the same CRP component, as we know the limiting service probabilities of customer classes and servers that are not in the same CRP component converge to zero. Similarly, we do not consider that service probabilities in any states $x$ not in $P_k(\sbb_k)$ for some $k \in [K']$ and $\sbb_k \in \Sigma_{\Sbb_k}$, as those states have idle servers, and hence have probabilities converging to zero. 

We will begin by writing the state dependent matching probability $\htp{q}_{ij}(x)$ for an arbitrary state $x \in P_{k(j)}(\sbb_{k(j)})$. We will let $j(x)$ denote the position in the server permutation of server $j$ in the state $x$ and similarly will let $j(s)$ denote the position of server $j$ in the server permutation $s$. We can look at $\htp{q}_{ij}(x)$ by conditioning on the position in the queuing network of the potential customer of type $i$ that $j$ serves. This lets us express $\htp{q}_{ij(x)}$ as 
\begin{align}
\htp{q}_{ij}(x) & = \sum_{r = j(x)}^{m} \left( \prod_{u = j(x)}^{r - 1} \frac{\lambda^{n_u}_{\{U(s_1,\dots,s_u \cap \overline{C(j)}\}}}{\lambda^{n_u}_{U(s_1,\dots,s_u)}} \right) \left(\lambda_i \sum_{y = 1}^{n_r} \frac{\lambda^{n_r-1}_{\{U(s_1,\dots,s_r) \cap \overline{C(j)}\}}}{\lambda^{n_r}_{U(s_1,\dots,s_r)}} \right)\nonumber \\
& = \lambda_i \sum_{r = j(x)}^{m} \left( \prod_{u = j(x)}^{r - 1} \frac{\lambda^{n_u}_{\{U(s_1,\dots,s_u) \cap \overline{C(j)}\}}}{\lambda^{n_u}_{U(s_1,\dots,s_u)}} \right) \left( \frac{\lambda^{n_r}_{U(s_1,\dots,s_r)} -\lambda^{n_r}_{\{U(s_1,\dots,s_r) \cap \overline{C(j)}\}} }{ \lambda^{n_r}_{U(s_1,\dots,s_r)} \left( \lambda_{U(s_1,\dots,s_r)} -\lambda_{\{U(s_1,\dots,s_r) \cap \overline{C(j)}\}}\right)} \right).
\end{align}

It will be useful to decompose this expression into two parts, $q^+_{ij}(x)$, the part of the expression representing a transition within the CRP component, and $q^0_{ij}(x)$, the part of the expression representing a transition outside of the CRP component. We suppress the dependence on $\epsilon$ to reduce clutter in the notation. So
\begin{equation*}
    q^+_{ij}(x) = \lambda_i \sum_{r = j(x)}^{m_k} \left( \prod_{u = j(x)}^{r - 1} \frac{\lambda^{n_u}_{\{U(s_1,\dots,s_u) \cap \overline{C(j)}\}}}{\lambda^{n_u}_{U(s_1,\dots,s_u)}} \right) \left( \frac{\lambda^{n_r}_{U(s_1,\dots,s_r)} -\lambda^{n_r}_{\{U(s_1,\dots,s_r) \cap \overline{C(j)}\}} }{ \lambda^{n_r}_{U(s_1,\dots,s_r)} \left( \lambda_{U(s_1,\dots,s_r)} -\lambda_{\{U(s_1,\dots,s_r) \cap \overline{C(j)}\}}\right)} \right),
\end{equation*}
and $q^0_{ij}(x) = \htp{q}_{ij}(x) - q^+_{ij}(x)$. Recall that $m_\kappa = \sum_{\ell \in [\kappa]} |\Scal_\ell|$, that is, $m_\kappa$ is the number of servers in the first $\kappa$ CRP components in the topological order. 

As an intermediate step to looking at the aggregate matching probabilities $\htp{q}_{ij}(P_k(\sbb_{k}))$, we will first look at the partially aggregated matching probabilities $\htp{q}_{ij}(P(s,m))$. 
\begin{align*}
    \htp{q}_{ij}(P(s,m)) & = \frac{1}{\pi(P(s,m))}\left[\sum_{x \in P(s,m)} \pi(x)q^+_{ij}(x) + \sum_{x \in P(s,m)} \pi(x)q^0_{ij}(x)\right].
\end{align*}
However, the second term represents transitions from a state where the permutation of servers is induced by a topological order to a state where the permutation of servers is not induced by a topological order, and hence has a limiting probability of zero. This means we expect the second term in this expression to converge to zero, which we prove in the following lemma.
\begin{lem}\label{lem:prob_transition_out}
For a given admissible service menu $M$ with limiting arrival rates $\Lambda$, service rates $\mu$, and slacks $\Gamma$, let $\{\mathbb{C}_1,\dots,\mathbb{C}_{K'}, \mathbb{C}_{K'+1}, \dots, \mathbb{C}_{K} \}$ be the set of CRP components, and let $\Tcal(\Dcal, K')$ be the set of topological orders on the CRP components. Then for any permutation of servers $s$ induced by some topological order $\sigma \in \Tcal(\Dcal, K')$, 
\begin{equation*}
    \lim_{\epsilon \to 0}\sum_{x \in P(s,m)} \pi(x)q^0_{ij}(x) = 0
\end{equation*}
\end{lem}
{\sc Proof}: See Appendix \ref{appn:proof_matching}. $\Box$

We will now fix a topological order $\sigma \in \Tcal(\Dcal, K')$, and a server permutation $s \in \Sigma_m$ that is induced by $\sigma$. To reduce notational clutter, we assume without loss of generality that the CRP components are labelled in order of their position in the topological order, that is, $\sigma(k) = k$ for all $k \in K'$. Using \cref{lem:prob_transition_out}, we can write $\htp{q}_{ij}(P(s,m))$ as
\begin{align*}
    \htp{q}_{ij}(P(s,m)) & = \frac{1}{\pi(P(s,m))}\sum_{x \in P(s,m)} \pi(x)q^+_{ij}(x) + o(1),
\end{align*}
or written another way,
\begin{align}\label{eq:matchingprob_permutation}
    \htp{q}_{ij}(P(s,m)) & = \frac{\lambda_i }{\pi(P(s,m))}\sum_{n_1 = 0}^{\infty} \cdots \sum_{n_m = 0}^ \infty {\cal B}  \prod_{\ell=1}^m \frac{\lambda^{n_\ell}_{U(s_1,\ldots, s_\ell)} }{\mu_{\{s_1,\ldots, s_\ell\}}^{n_\ell+1} } q^+_{ij}(s_1, n_1, \cdots, s_m, n_m) + o(1).
\end{align}

The following notation will be useful in simplifying this expression. Recall from \cref{def:capacity_slack} that 
\begin{align*}
    \Delta(S) &= \mu_S - \lambda_{U_S(M)}. 
\end{align*}
It will also be useful to define $\Delta_j(S)$ as 
\begin{align}
    \Delta_j(S) &= \mu_S - \lambda_{\{U_S(M) \cap \overline{C(j)}\}}. 
\end{align}

We can then write \cref{eq:matchingprob_permutation} as
\begin{align}\label{eq:matchingprob_permutation_2}
\htp{q}_{ij}(P(s,m)) & = \frac{\mathcal{B}\lambda_i}{\pi(P(s,m))} \left(\prod_{\ell = m_{k(j)}+1}^{m}\frac{1}{\Delta(s_1, \ldots, s_\ell)} \right) \left(\prod_{\ell = 1}^{m_{k(j)} - 1}\frac{1}{\Delta(s_1, \ldots, s_\ell)} \right) \nonumber \\
 & \times \left(\prod_{\ell = m_{k(j)-1} + 1}^{j-1} \frac{1}{\Delta(s_1, \ldots, s_\ell)} \right) \Bigg[   \sum_{r = j(s)}^{m_k(j)} \left(\prod_{u = j(s) }^{r} \frac{1}{\Delta_j(s_1, \ldots, s_u)} \right) \left(\prod_{\ell = r + 1 }^{m_{k(j)}} \frac{1}{\Delta(s_1, \ldots, s_\ell)} \right) \nonumber \\
& \times \left( \frac{1}{\Delta(s_1, \ldots, s_r)} - \frac{1}{\Delta_j(s_1, \ldots, s_r)} \right) \Bigg] + o(1),
\end{align}
where as before $m_\kappa = \sum_{\ell \in [\kappa]} |\Scal_\ell|$.That is, $m_\kappa$ is the number of servers in the first $\kappa$ CRP components in the topological order. 

We saw in \cref{sec:waits_proof} that the limiting values of $\Delta(s_1, \ldots, s_\ell)$ depend on the values of $\ell$. If $\ell = m_\kappa$ for some $\kappa \in [K']$, then we know from \cref{lem:criticalservers} that 
\[ \lim_{\epsilon \to 0} {\epsilon \over \Delta(s_1, \ldots, s_{m_\kappa})} = \frac{1}{ \sum_{\ell = 1}^\kappa \widetilde{\gamma}_{\comps(\sigma, \ell)}}. \]

For all other values of $\ell$, there is some $\kappa \in [K']$ such that $ m_{\kappa-1} + 1 \leq \ell \leq m_{\kappa}-1$. Here we take $m_0 = 0$. We let $S = \{ s_{m_{\kappa-1}+1}, \ldots, s_\ell\}$. Following the outline in \cite[Lemmas 5 and 8]{Afecheetal2019}, we can show that:
\begin{align*}
    \lim_{\epsilon \to 0} \Delta(s_1,\ldots, s_{\ell}) &= \mu_{S} - \Lambda_{U_{S}(\breve{M}) } > 0.
\end{align*} 
In particular, this means that for all $\kappa \in [K']$, and $ m_{\kappa-1} + 1 \leq \ell \leq m_{\kappa}-1$,  $\lim_{\epsilon \to 0} \Delta(s_1,\ldots, s_{\ell})$ is a real number greater than zero that depends only on the permutation of servers in $\mathbb{C}_\kappa$. 

The same reasoning implies that for all $j \le \ell \le m_{k(j)}$, $\lim_{\epsilon \to 0} \Delta_j(s_1,\ldots, s_{\ell})$ is a real number greater than zero that depends only on the permutation or servers in  $\mathbb{C}_k$. 

We can use these observations to prove the following lemma.
\begin{lem}\label{lem:q_ij_limit} 
We can find functions $\left\{ \theta_\kappa : \Sigma_{\Scal_\kappa} \to \Re^+ \right\}_{\kappa \in [K']}$, $H_{ij} : \Sigma_{\Scal_{k(j)}} \to \Re^+$, and $ G_{ij} : \Sigma_{\Scal_{k(j)}} \to \Re^+$, such that $q_{ij}(P(s,m)) = \lim_{\epsilon \to 0}\htp{q}_{ij}(P(s,m))$ can be written as 
\begin{align}\label{eq:q_ij_limit}
q_{ij}(P(s,m)) = & \lim_{\epsilon \to 0} \left[ \frac{\Bcal\lambda_i}{\pi(P(s,m))\epsilon^{K'}} \mathbb{Q}(\sigma) \left(\prod_{\kappa \neq k(j)} \theta_\kappa(s_\kappa) \right) H_{ij}(\sbb_{k(j)}) \right] \nonumber \\
& - \lim_{\epsilon \to 0} \left[ \frac{\Bcal\lambda_i}{\pi(P(s,m))\epsilon^{K'-1}} \left(\prod_{\kappa \neq k} \frac{1}{\sum_{\ell = 1}^\kappa \widetilde{\gamma}_{\comps(\sigma, \ell)}} \right) \left(\prod_{\kappa \neq k(j)} \theta_\kappa(s_\kappa) \right) G_{ij}(\sbb_{k(j)}) + o(1)\right],
\end{align}
where $\theta_\kappa$ and $H_{ij}$ only depend on $\breve{M}$, $\Lambda$, and $\mu$. 
\end{lem}
{\sc Proof}: See Appendix \ref{appn:proof_matching}. $\Box$

We provide exact definitions of $\left\{ \theta_k : \Sigma_{\Scal_k} \to \Re^+ \right\}_{k \in [K']}$, $H_{ij} : \Sigma_{\Scal_k} \to \Re^+$, and $ G_{ij} : \Sigma_{\Scal_k} \to \Re^+$ in the proof of \cref{lem:q_ij_limit} in Appendix \ref{appn:proof_matching}.

Notice that the first line in \cref{eq:q_ij_limit} has an $\epsilon^{-K'}$ term, and the second line has an $\epsilon^{-(K'-1)}$ term. Since $q_{ij}$ are probabilities and therefore must be between 0 and 1, we know that $\lim_{\epsilon \to 0}B\epsilon^{-K'}$ is bounded. This implies that $\lim_{\epsilon \to 0}B\epsilon^{-(K'-1)} = 0$, and so only the first line in \cref{eq:q_ij_limit} will be non-zero. Thus 
\begin{align}
q_{ij}(P(s,m)) & = \frac{\Bcal'\lambda_i}{\pi(P(s,m))} \mathbb{Q}(\sigma) \left(\prod_{\kappa \neq k} \theta_\kappa(s_\kappa) \right) H_{ij}(\sbb_{k(j)}).
\end{align}

Because $q_{ij}$ are matching probabilities, we also know that 
\begin{align}\label{eq:q_ij_limit2}
    q_{ij}(P(s,m)) & = \frac{q_{ij}(P(s,m))}{\sum_{i' \in \Ccal_{k(j)}}q_{i'j}(P(s,m)) }.
\end{align}
Since the only term in \cref{eq:q_ij_limit2} that depend on $j$ is the $H_{ij}(\sbb_{k(j)})$ term, we can write $q_{ij}(P(s,m))$ as 
\begin{align}\label{eq:q_ij_limit_3}
q_{ij}(P(s,m)) = \frac{H_{ij}(\sbb_{k(j)})}{\sum_{i' \in \Ccal_k}H_{i'j}(\sbb_{k(j)})}. 
\end{align}
Since \cref{eq:q_ij_limit_3} holds for any server permutation $s \in \Sigma$, and depends only on $\sbb_k$ and not on the rest of the server permutation, this implies that 
\begin{align}\label{eq:q_ij_limit_4}
q_{ij}(P_{k(j)}(\sbb_{k(j)})) = \frac{H_{ij}(\sbb_{k(j)})}{\sum_{i' \in \Ccal_k}H_{i'j}(\sbb_{k(j)})}. 
\end{align}

Since, as \cref{lem:q_ij_limit} states, $H_{ij}(\sbb_{k(j)})$ does not depend on $\Gamma$, the remaining step needed to prove \cref{thm:matching_independent_gamma} is to show that $\pi(P_{k(j)}(\sbb_{k(j)}))$ also does not depend on $\Gamma$. This is captured in the following lemma. 
\begin{lem}\label{lem:permutation_prob}
For an admissible service menu $M$ with limiting arrival rates $\Lambda$ service rates $\mu$, and slacks $\Gamma$, the limiting probability of being in a state with the sub-permutation of server $\sbb_k \in \Sigma_{\Sbb_k}$ for $k \in K'$ is equal to 
\begin{equation*}
    \lim_{\epsilon \to 0}\pi(P_k(\sbb_k)) =  \frac{\theta_k(s_k)}{ \sum_{\sbb_\kappa \in \Sigma_{\Scal_k} } \theta_\kappa( \sbb_\kappa)},
\end{equation*}
where $\left\{ \theta_\kappa : \Sigma_{\Scal_\kappa} \to \Re^+ \right\}_{\kappa \in [K']}$ is a function that depends only on $\breve{M}$, $\Lambda$, and $\mu$.
\end{lem}
{\sc Proof}: See Appendix \ref{appn:proof_matching}. $\Box$

Combining \cref{lem:permutation_prob} with \cref{eq:q_ij_limit_3}, we have that the limiting service probabilities $\lim_{\epsilon \to 0}\htp{q}_{ij}$ do not depend on the exact values of the slacks $\Gamma$, only requiring that $M$ is an admissible menu for the slacks $\Gamma$.

\section{Concluding Remarks}\label{sec:Conlusion}
In this paper, we have studied the performance of multi-class multi-server bipartite queueing systems under a FCFS-ALIS service discipline by extending the heavy traffic analysis introduced in \cite{Afecheetal2019} for a similar class of systems. In \cref{thm:CRPdelay_append} we have provided a general characterization of the mean steady-state waiting time delay for each customer class. Our characterization relies on decomposing the queueing system into a collection of complete resource pooling (CRP) components and identifying the connectivity among these CRP components in the form of a directed acyclic graph (DAG). Interestingly, only the knowledge of this DAG together with the capacity slack in each CRP component is enough to derive the mean steady-state waiting time for all customer classes. We have also studied the steady-state matching probabilities among customer classes and servers and showed in \cref{thm:matching_independent_gamma} that only the limiting values of arrival and service rates influence these matching probabilities. This is in direct contrast to the behaviour of the mean steady-state waiting times, which are also affected by the direction of convergence to heavy traffic. To illustrate this point, we have provided a numerical example that shows that small changes to the arrival rates in a heavily congested system can have large impacts on the average delays.  We use our results regarding steady-state outcomes to explore some questions regarding the design of queueing systems. In doing this, we find that when service providers are looking to minimise expected delays and have complete control over the design of the menu, then they should implement a menu that induces a single CRP component. 

Our work points towards several promising research directions. Firstly, we suggest exploring the 
problem of menu design, which involves determining the service classes to offer when customers 
can select which queue to join upon arrival. \cite{caldenteyetal2023} have made some preliminary progress in this area. Another area that deserves further investigation is the relationship 
between delays and the underlying matching topology in our bipartite queueing system. In \cref{sec:numex}, we demonstrate that adding more connectivity to the system can lead to a deterioration in the average waiting time of customers, exhibiting a form of Braess's paradox, despite neither customers nor servers acting strategically. 
Mathematically, this negative effect happens when adding an additional arc to the menu increases the probability of a topological order with higher conditional delays. \cref{thm:CRPdelay_append} characterizes waiting time delays and can be used to identify an optimal flexibility structure as a combinatorial optimization problem over the collection of directed acyclic graphs (DAGs) associated with a particular set of CRP components.

In addition, there are alternative modelling choices that could be worth exploring. For example, while we have focused on conventional heavy-traffic scaling in this paper, a many-server scaling may be more appropriate for certain application settings, such as public housing and healthcare, where many identical servers are available. Furthermore, we have primarily examined steady-state outcomes, but in real-world scenarios, conditions often change frequently, making it unclear if a steady-state will be achieved. Therefore, studying the transient behaviour of bipartite queueing systems could also be of interest.

\bibliographystyle{abbrvnat}
\bibliography{BiblioMatching.bib}

	\clearpage
\begin{appendices}
 
\counterwithin{equation}{section}
\renewcommand{\theequation}{A\arabic{equation}}
\setcounter{equation}{0}
\renewcommand\thesection{\Alph{section}:}

\section{\cref{sec:heavytrafficwaits} Proofs} \label{sec:proof_W}
{\sc Proof of Lemma~\ref{lem:flows}:} Let us define the set $\mathcal{F}_{\max}$ as 
\begin{equation*}
    \mathcal{F}_{\max} := \left\{ \sum_{i \in [n]}f = [f_{ij}] : \quad \sum_{i \in [n]}f_{ij} \le \mu_j \quad \forall j \in [m] \quad , f \ge 0, \quad f_{ij}= 0, ~~\forall (i,j) : m_{ij}=0  \right\}.
\end{equation*}
Note that for all $\epsilon \in [0, \epsilon_o)$, 
$\mathcal{F}(\epsilon, \htp{\lambda}, M) \subseteq \mathcal{F}_{\max}$. Furthermore, since $\mathcal{F}_{\max}$ is a compact set, we know that the sequence $\htp{f}$ 
has a subsequence that converges to some limit in $\mathcal{F}_{\max}$. Let $\tilde{f}$ denote this limit. To prove that $\tilde{f} \in \mathcal{F}(0, \htp{\lambda}, M)$, all that remains to be shown is that $\tilde{f}$ satisfies
\begin{equation*}
      \sum_{j \in [m]}\tilde{f}_{ij} = \Lambda_i, \quad \text{for all } i \in [n].
\end{equation*}
But we know that 
\begin{equation*}
    \sum_{j \in [m]}\htp{f}_{ij} = \htp{\lambda}_i, \quad \text{for all } i \in [n] \text{ and } 0 < \epsilon < \epsilon_0,
\end{equation*}
and $\tilde{f}$ is the limit of a subsequence of $\htp{f}$, and so 
\begin{equation*}
    \sum_{j \in [m]}\tilde{f}_{ij} = \lim_{\epsilon \to 0} \htp{\lambda}_i = \Lambda_i, \quad \text{for all } i \in [n]
\end{equation*}
as required. \hfill $\Box$
\ \\

{\sc Proof of Lemma~\ref{lem:slack_positive}:} Fix a topological order $\sigma_t \in \Tcal(\Dcal, [K'])$ and an index $\kappa \in [K']$. Define the sets
\[ \mathscr{C}= \bigcup\limits_{\ell=1}^{\kappa} \left\{ \Ccal_i : i \in  \comps(\sigma_t,\ell)  \right\}, \quad \mbox{and}\quad \mathscr{S} = \bigcup\limits_{\ell=1}^{\kappa} \left\{ \Scal_i : i \in  \comps(\sigma_t,\ell)  \right\}.\]
By the definition of the DAG $\Dcal$ and topological order $\sigma_t$, we have that
\[  \mathscr{S} = S(\mathscr{C}) .\]
That is, the services classes $\mathscr{C}$ are only served by servers in $\mathscr{S}$. We can find a lower bound on the scaled mean waiting times of the service classes in $\mathscr{C}$ using the scaled mean waiting time of a $M/M/1$ queue:
\begin{align}
\label{eq:mean_W_CRP_prefix}
    \sum_{i \in \mathscr{C}} \htp{\lambda}_i \hts{W_i} \geq \frac{\epsilon}{\mu_{\mathscr{S}} - \htp{\lambda}_{\mathscr{C}}}. 
\end{align}
Further, from \cref{lem:CRP_properties} we know that, 
\[  \mu_{\mathscr{S}} - \htp{\lambda}_{\mathscr{C}} = \epsilon \sum_{\ell = 1}^\kappa \widetilde{\gamma}_{\comps(\sigma_t, \ell)} + o(\epsilon).\]
If, contradictory to the \cref{lem:slack_positive}, $\sum_{\ell = 1}^\kappa \widetilde{\gamma}_{\comps(\sigma_t, \ell)} \leq 0$, then the right-hand side of \eqref{eq:mean_W_CRP_prefix} must diverge, and hence the sum on the left-hand side as well. However, from the admissibility of $M$, each $\hts{W_i}$ converges, and therefore also the sum on the left-hand side of \eqref{eq:mean_W_CRP_prefix}. Thus we must have $\sum_{\ell = 1}^\kappa \widetilde{\gamma}_{\comps(\sigma_t, \ell)} > 0$ for all $\sigma_t \in \Tcal(\Dcal, [K'])$ and $\kappa \in [K']$. \hfill $\Box$
\ \\
\section{\cref{sec:matchingprob} Proofs} \label{sec:proof_matchingprob_statement}
{\sc Proof of Proposition~\ref{prop:limiting_admissibility}}: Given $M$ is admissible for $(\htp{\lambda}, \mu)$, we know from the definition of admissibility that 
\begin{equation}\label{eq:admiss_proof_given}
\htp{\Delta}_{\mathscr{S}}(M) := \sum_{j \in \mathscr{S}} \mu_j - \sum_{i \in U_{\mathscr{S}}(M)} \htp{\lambda}_i = \Omega(\epsilon) \qquad \mbox{for all } \mathscr{S} \subseteq [m].
\end{equation}
To show $M$ is admissible for $(\Lambda - \epsilon\Lambda, \mu)$, we must show that 
\begin{equation}\label{eq:admiss_proof_toshow}
\sum_{j \in \mathscr{S}} \mu_j - \sum_{i \in U_{\mathscr{S}}(M)} \Lambda_i + \epsilon \sum_{i \in U_{\mathscr{S}}(M)} \Lambda_i= \Omega(\epsilon) \qquad \mbox{for all } \mathscr{S} \subseteq [m].
\end{equation}
\cref{eq:admiss_proof_given} implies that $\sum_{j \in \mathscr{S}} \mu_j - \sum_{i \in U_{\mathscr{S}}(M)} {\Lambda}_i \ge 0$ for all $\mathscr{S} \subseteq [m]$. For any $\mathscr{S} \subseteq [m]$ such that $\sum_{j \in \mathscr{S}} \mu_j - \sum_{i \in U_{\mathscr{S}}(M)} {\Lambda}_i > 0$, \cref{eq:admiss_proof_toshow} holds without regardless of the $\epsilon$ terms. In the case that $\sum_{j \in \mathscr{S}} \mu_j - \sum_{i \in U_{\mathscr{S}}(M)} {\Lambda}_i = 0$, then 
\begin{align*}
\sum_{j \in \mathscr{S}} \mu_j - \sum_{i \in U_{\mathscr{S}}(M)} \Lambda_i + \epsilon \sum_{i \in U_{\mathscr{S}}(M)} \Lambda_i & = \epsilon \sum_{i \in U_{\mathscr{S}}(M)} \Lambda_i \\
& = \epsilon \sum_{j \in \mathscr{S}} \mu_j. 
\end{align*}
But $\sum_{j \in \mathscr{S}} \mu_j > 0$, so $\epsilon \sum_{j \in \mathscr{S}} \mu_j= \Omega(\epsilon)$ as required. 

The second part of the proposition states that $\breve{M}$ is admissible for $(\Lambda - \epsilon \Lambda, \mu)$. To show this, similarly to the first part of the proposition we must show that 
\begin{equation}\label{eq:proof_admiss_residual}
\sum_{j \in \mathscr{S}} \mu_j - \sum_{i \in U_{\mathscr{S}}(\breve{M})} \Lambda_i + \epsilon \sum_{i \in U_{\mathscr{S}}(\breve{M})} \Lambda_i= \Omega(\epsilon) \qquad \mbox{for all } \mathscr{S} \subseteq [m].
\end{equation}
There are two cases to consider. In the first case, $\mathscr{S} = \cup_{k \in T} \mathcal{S}_k$ for some $T \subseteq K$. In words, this means that $\mathscr{S}$ is the union of servers in a particular subset of CRP components. It is shown in the proof of \cite[Lemma 4]{Afecheetal2019} that in this case, $\sum_{j \in \mathscr{S}} \mu_j = \sum_{i \in U_{\mathscr{S}}(\breve{M})} \Lambda_i$, and \cref{eq:proof_admiss_residual} holds following the same reasoning as in the first part of the proposition. In the second case, $\mathscr{S} \neq \cup_{k \in T} \mathcal{S}_k$ for any $T \subseteq K$, and the proof of \cite[Lemma 4]{Afecheetal2019} shows that $\sum_{j \in \mathscr{S}} \mu_j > \sum_{i \in U_{\mathscr{S}}(\breve{M})} \Lambda_i$, and \cref{eq:proof_admiss_residual} holds following similar reasoning as in the first part of the proposition. \hfill $\Box$

\section{\cref{sec:discussion} Proofs} \label{sec:proof_discussion}
{\sc Proof of Corollary~\ref{cor:implementable}}: We will prove this corollary by proving the contrapositive. So suppose there are $k \in [K]$ and $\kappa \in [K]$ such that there are no topological orders $\sigma \in \Tcal(\Dcal, K')$ with $\comps^{-1}(\sigma, \kappa) \le \comps^{-1}(\sigma, k)$. This means that in every topological order $\sigma \in \Tcal(\Dcal, K')$, $\comps^{-1}(\sigma, \kappa) > \comps^{-1}(\sigma, k)$. From the definition of the conditional delay $w_{\sigma, k}$ in \cref{eq:w_sigma_k}, this implies that $w_{\sigma, k} > w_{\sigma, \kappa}$ for all $\sigma \in \Tcal(\Dcal, K')$. As the total delays are weighted sums of the conditional delays, this proves the result. 
\hfill $\Box$
\ \\

{\sc Proof of Proposition~\ref{prop:Wchained}}: Without loss of generality let us index the CRP components in such a way that $W_k\le W_{k+1}$ for all $k \in [K - 1]$. Recall $\{\mathscr{C}_1, \dots, \mathscr{C}_L\}$ is the partition described in \cref{def:chain}. Are stated in the proposition, we will assume 
\begin{enumerate}[nosep]
    \item[\rm (i)] $W_k = W_\kappa$ for all $(k, \kappa) \in [K] \times [K]$ such that $W_k \in \mathscr{C}_\ell$ and $W_\kappa \in \mathscr{C}_\ell$ for some $\ell \in [L]$,
    \item[\rm (ii)]  $W_k < W_\kappa$ for all $(k, \kappa) \in [K] \times [K]$ such that $W_k \in \mathscr{C}_\ell$ and $W_\kappa \in \mathscr{C}_{\ell'}$ for some $(\ell, \ell') \in [L] \times [L]$ where $\ell < \ell'$.
\end{enumerate}
We will now show how to choose a vector of capacity slacks $\tilde{\gamma}=(\tilde{\gamma}_1,\dots, \tilde{\gamma}_K)$ such that $W_{\mathbb{C}_k} = W_k$ for all $k \in [K]$. Fix $\tilde{\gamma}$ such that $\widetilde{\gamma}_k=\widehat{\gamma}_\ell$ for all $k \in \mathscr{C}_\ell$. It follows from the chained structure of the DAG and the construction of $\tilde{\gamma}$ that for any permutation  $\sigma=(\sigma(1),\sigma(2),\dots,\sigma(K))$  induced by some topological order the vector $(\widetilde{\gamma}_{\sigma^{-1}(1)}, \widetilde{\gamma}_{\sigma^{-1}(2)}, \dots, \widetilde{\gamma}_{\sigma^{-1}(K)})$ is constant. This observation together with \cref{thm:CRPdelay_append} imply that $\mathbb{Q}(\sigma)$ in \cref{eq:Q_def} is also constant, independent of $\sigma$.  Furthermore, by symmetry, it is not hard to see that two CRP components that belong to the same partition $\mathscr{C}_\ell$ have the same limiting scaled waiting times, which we denote by $\widehat{\mathbb{W}}_\ell$. One can show from \cref{thm:CRPdelay_append} that 
\begin{equation}\label{eq:appenchained}\widehat{\mathbb{W}}_{\ell} =\widehat{\mathbb{W}}_{\ell-1}+{1 \over n_\ell} \sum_{s=1}^{n_\ell} {1 \over \sum_{j=\ell+1}^L n_j\, \widehat{\gamma}_j+s\,\widehat{\gamma}_\ell}, \qquad \ell=1,2\dots,L\end{equation}
with $\widehat{\mathbb{W}}_{0}=0$. We use this condition to find the values of $\{\widehat{\gamma}_\ell\}$ that implement $\{\mathbb{W}_\ell\}$, that is,  $\widehat{\mathbb{W}}_\ell=\mathbb{W}_\ell$ for all $\ell \in [L]$. To this end, we use backward induction on $\ell$. For $\ell=L$ we have that
$$\widehat{\mathbb{W}}_L=\widehat{\mathbb{W}}_{L-1}+{1 \over n_L} \sum_{s=1}^{n_L} {1 \over s\,\widehat{\gamma}_L}.$$ Thus, $\widehat{\gamma}_L$ must satisfy
$$\widehat{\gamma}_L = {1 \over (\mathbb{W}_L-\mathbb{W}_{L-1})}\, {1 \over n_L} \sum_{s=1}^{n_L} {1 \over s}.$$
Now suppose that we have determined the values of $\widehat{\gamma}_L, \widehat{\gamma}_{L-1}, \dots, \widehat{\gamma}_{\ell+1}$ and define $\widehat{\Gamma}_\ell:=\sum_{j=\ell+1}^L n_j\,\widehat{\gamma}_j.$  We find the value $\widehat{\gamma}_{\ell}$ by solving \eqref{eq:appenchained}
$$\mathbb{W}_{\ell} =\mathbb{W}_{\ell-1}+{1 \over n_\ell} \sum_{s=1}^{n_\ell} {1 \over \widehat{\Gamma}_\ell+s\,\widehat{\gamma}_\ell}.$$ We note that there exists a unique $\widehat{\gamma}_{\ell}$ that solves this equation in the region $\widehat{\gamma}_{\ell} >-\widehat{\Gamma}_\ell / n_\ell$. This follows from the fact that the summation above is monotonically decreasing in $\widehat{\gamma}_{\ell}$ in this region and diverges to $+\infty$ as $\widehat{\gamma}_{\ell}$ approaches $\-\widehat{\Gamma}_\ell / n_\ell$ from above and converges to zero as $\widehat{\gamma}_{\ell}$ approaches $\infty$. 
\hfill $\Box$
\ \\

{\sc Proof of Proposition~\ref{prop:min_avg_waits}}: Note from \eqref{eq:Q_def} that 
$$ w_{\sigma, k} := \sum_{ \kappa =  \sigma^{-1}(k)  }^{K} \frac{1}{\sum_{\ell=1}^\kappa \widetilde{\gamma}_{\sigma(\ell)}} = {1 \over |a|}+ \sum_{ \kappa =  \sigma^{-1}(k)  }^{K-1} \frac{1}{\sum_{\ell=1}^\kappa \widetilde{\gamma}_{\sigma(\ell)}}.$$
Let us prove that $ w_{\sigma, k} \geq 1/|a|$. From the previous equation, this would follow if the last summation is nonnegative. Suppose, by contradiction that this is not the case. Then, there exists a $\kappa$ such that $\sigma^{-1}(k)  \leq \kappa \leq K-1$ such that $\sum_{\ell=1}^\kappa \widetilde{\gamma}_{\sigma(\ell)}<0$. In other words, the cumulative capacity slack of the CRP components $\{\mathbb{C}_{\sigma(1)},\mathbb{C}_{\sigma(2)},\dots,\mathbb{C}_{\sigma(\kappa)}\}$ is negative. However, this would imply that the cumulative arrival rate to these components exceeds the total service capacity of all the servers in these components. This, together with the DAG structure connecting all the CRP components imply that the stability condition in \cref{prop:stability} is violated, which holds by assumption. From this contradiction we conclude that $ w_{\sigma, k} \geq 1/|a|$ and then from \eqref{eq:WCRP_def} we also get that $W_{\mathbb{C}_k}  \geq 1/|a|$.

Let us now prove the second part of the corollary, namely, there can be at most one CRP component $\hat{\kappa} \in [K]$ such that $\widehat{W}_{\mathbb{C}_{\hat{\kappa}}}=1/|a|$.  From the previous discussion, it follows that the requirement $\widehat{W}_{\mathbb{C}_{\hat{\kappa}}}=1/|a|$ can only be satisfied if $ w_{\sigma, \hat{\kappa}}= 1/|a|$ for all permutations $\sigma$ associated a topological order.   But this can only happen if  $\sigma^{-1}(\hat{\kappa})=K$ for all permutation $\sigma$. Evidently, this condition can only be satisfied by at most one CRP component and holds trivially if $K=1$. $\Box$

{\sc Proof of Proposition~\ref{prop:singleCRP_condition}}: Take any slacks $\gamma$ with $|\gamma| > 0$. We will first show that $M$ is admissible with $\htp{\lambda} = \Lambda - \epsilon\gamma + o(\epsilon)$ and $\mu$. To do this, we need to show that 
\begin{equation*}
    \htp{\Delta}_{\mathscr{S}}(M) = \Omega(\epsilon) \qquad \mbox{for all } \mathscr{S} \subseteq [m], 
\end{equation*}
where 
\begin{equation*}
\htp{\Delta}_{\mathscr{S}}(M) := \sum_{j \in \mathscr{S}} \mu_j - \sum_{i \in U_{\mathscr{S}}(M)} \htp{\lambda}_i.
\end{equation*}
We define $D_{\mathscr{S}}$ as 
\begin{equation*}
    D_{\mathscr{S}} = \sum_{j \in \mathscr{S}} \mu_j - \sum_{i \in U_{\mathscr{S}}(M)} \Lambda_i
\end{equation*}
for all $\mathscr{S} \subseteq [m]$. Then 
\begin{equation*}
\htp{\Delta}_{\mathscr{S}}(M) = D_{\mathscr{S}} + \epsilon\sum_{j \in \mathscr{S}}\gamma_i + o(\epsilon) \qquad \mbox{for all } \mathscr{S} \subseteq [m], 
\end{equation*}
From the definition of $M$ we know that $D_{\mathscr{S}} > 0$ for all $\mathscr{S} \subseteq [m]$, implying that $\htp{\Delta}_{\mathscr{S}}(M) = \Omega{\epsilon}$ for all $\mathscr{S} \subseteq [m]$. For the case of $\mathscr{S} = [m]$, since $|\Lambda| = |\mu|$, and $|\gamma| > 0$, 
\begin{equation*}
\htp{\Delta}_{\mathscr{S}}(M) = +\epsilon|\gamma| + o(\epsilon) = \Omega(\epsilon)
\end{equation*}
as required. 

What remains to be shown is that $M$ induces a single CRP component. This follows from part (i) of \cref{lem:CRP_properties}, which states that within a CRP component $\widetilde{\Lambda}_k := \Lambda_{\Ccal_k} = \mu_{\Scal_k} =: \widetilde{\mu}_k$ (see \eqref{eq:CRP_mu_lambda} for definitions). But with our choice of $M$, we know that for any subset of servers $\mathscr{S} \subsetneq [m]$, any subset of customers classes $\mathscr{C} \subseteq [n]$ such that every class in $\mathscr{C}$ is compatible with some server in $\mathscr{S}$ will have $\Lambda_{\mathscr{C}} < \mu_{\mathscr{S}}$. Thus there are no CRP components that do not consist of all customer classes and all servers, implying there is exactly one CRP component. 

\hfill $\Box$

{\sc Proof of Lemma~\ref{lem:admissible_top}}: We assume without loss of generality that the CRP components are labelled so that $\comps^{-1}(\sigma, k) = k$ for all $k \in [K']$. We construct the menu $M$ as follows. Let $\breve{M}$ be any residual matching associated with the collection of CRP components $\mathbb{C} = \{\mathbb{C}_1,\dots,\mathbb{C}_{K'},\mathbb{C}_{K'+1},\dots,\mathbb{C}_{K}\}$. Construct the menu $M$ as follows. Let $m_{ij} = 1$ for all $i \in [n]$ and $j \in [m]$ such that $\breve{m}_{ij} = 1$. Then for every $k \in [K' - 1]$, let $m_{ij} = 1$ for some $i \in \mathcal{C}_{k+1}$ and some $j \in \mathcal{S}_k$. That is, for every CRP component $\mathbb{C}_k$ for $k \in [K'-1]$, we assign some customer class in $\mathbb{C}_{k+1}$ to be a served by a server in $\mathbb{C}_{k}$. We will show that this has the effect of adding an arc to the DAG from $\mathbb{C}_{k+1}$ to $\mathbb{C}$ without altering the CRP component structure. 

We will begin by assuming that there are no customer classes with zero arrivals, that is, we assume that $\tilde{\Lambda}_k > 0$ for all $k \in [K]$, and $K = K'$. In this case, we let $m_{ij} = 0$ for all other combinations of $i \in [n]$ and $j \in [m]$. We will mention at the end of this proof how to adjust the menu $M$ for the case in which there is at least one $k \in [K]$ with $\tilde{\Lambda}_k = 0$. 

The next step is to show that the CRP components of $M$ are $\mathbb{C}$. This is equivalent to showing that $\mathcal{F}(0, \Lambda, M) = \mathcal{F}(0, \Lambda, \breve{M})$. First note that there can only be flow between servers in $\mathcal{S}_k$ and customers in $\mathcal{C}_{k} \cup \mathcal{C}_{k+1}$ for $k \in [K - 1]$, and there can only be flow between servers in $\mathcal{S}_{K}$ and customers in $\mathcal{C}_K$ due to the construction of $M$. But there can be no flow between servers in $\mathcal{S}_1$ and customers in $\mathcal{C}_2$, as all of the capacity of servers in $\mathcal{S}_1$ needs to be allocated to servers in $\mathcal{C}_1$, since $\tilde{\Lambda}_1 = \tilde{\mu}_1$. It can then be argued inductively that servers in $\mathcal{S}_k$ do not have the capacity to allocate flow to customers in $\mathcal{C}_{k+1}$, even though there is a server that has the compatibility to do so. Thus $\mathcal{F}(0, \Lambda, M) = \mathcal{F}(0, \Lambda, \breve{M})$ as required. 

Next, we will show that the DAG of $M$ only admits the topological order $\sigma$. This is true based on the construction of $M$. The only arcs in $M$ that are not in the residual matching $\breve{M}$ are between components $\mathbb{C}_k$ and $\mathbb{C}_{k+1}$ for $k \in [K - 1']$, and there is such an arc for $k \in [K - 1]$. Thus we require for any topological order $\sigma_t$ admitted by $M$ that $\sigma_t(k) < \sigma_t(k+1)$ for $k \in [K - 1]$. But the only topological order that achieves this is $\sigma$, where as stated previously $\sigma(k) = k$. 

The final step needed to prove the first claim in \cref{lem:admissible_top} is to show that $M$ is admissible. Recall from \cref{def:admissible} that for a menu to be admissible we require that $\htp{\Delta}_{\mathscr{S}}(M) = \Omega(\epsilon)$ for all $\mathcal{S} \subseteq[m]$, where 
\begin{equation*}
    \htp{\Delta}_{\mathscr{S}}(M) := \sum_{j \in \mathscr{S}} \mu_j - \sum_{i \in U_{\mathscr{S}}(M)} \htp{\lambda}_i.
\end{equation*}
The proof of \cite[Lemma 4]{Afecheetal2019} argues that if the subset of servers $\mathscr{S} \subseteq [m]$ is not equal to $\cup_{\ell = 1}^k \mathcal{S}_\ell$ for some $k \in [K]$, then
\[   \mu_S - \Lambda_{U_S(M)} > 0. \]
which means that $\htp{\Delta}_{\mathscr{S}}(M) = \Omega(\epsilon)$ for all $\mathscr{S} \subseteq [m]$ that is not equal to $\cup_{\ell = 1}^k \mathcal{S}_\ell$ for some $k \in [K]$. For $\mathscr{S} \subseteq [m]$ such that $\mathscr{S} = \cup_{\kappa = 1}^k \mathcal{S}_\ell$ for some $k \in [K]$, we know from \cref{lem:CRP_properties} that $\sum_{j \in \mathscr{S}} \mu_j = \sum_{i \in U_{\mathscr{S}}(M)} \Lambda_i$. So 
\begin{equation*}
    \htp{\Delta}_{\mathscr{S}}(M) \sum_{\ell = 1}^k\epsilon\tilde{\gamma}_\ell - o(\epsilon). 
\end{equation*}
But since from the statement of the lemma, $\sum_{\ell = 1}^k\epsilon\tilde{\gamma}_\ell > 0$ for all $k \in [K]$, this means that $\htp{\Delta}_{\mathscr{S}}(M) = \Omega(\epsilon)$ as required. Hence $M$ is admissible as claimed.

This also demonstrates why no admissible menu $M$ can admit a topological order $\sigma$ such that $\sum_{\ell = 1}^k\epsilon\tilde{\gamma}_{\ell} \le 0$ for some $k \in [K']$. If that were the case, then we would have that $\lim_{\epsilon \to 0}\htp{\Delta}_{\mathscr{S}}(M) \le 0$ for $\mathscr{S} =\cup_{\kappa = 1}^k \mathcal{S}_\ell$, which contradicts $M$ being admissible. This holds even if we were to consider the scenario in which $\tilde{\Lambda}_k = 0$ for some $k \in [K]$, as this would only decrease the values of $\tilde{\gamma}_{\comps(\sigma), k}$, making it more difficult to satisfy the condition $\lim_{\epsilon \to 0}\htp{\Delta}_{\mathscr{S}}(M) > 0$. 

Finally, we will mention how we can extend the construction of $M$ to account for CRP components $k$ with $\tilde{\Lambda}_k = 0$. Recall that these CRP components do not influence the topological orders themselves, only the slacks the elements $\comps(\sigma, k)$. We require for the admissibility of $M$ that $\sum_{\ell = 1}^k\comps(\sigma, \ell) > 0$ for all $k \in [K']$. This can potentially be achieved in many ways, one of which will always be to let $m_{ij} = 1$ for some $j$ in $\mathbb{C}_{K'}$ and for all $i \in [n]$ such that $\Lambda_{i} = 0$. This construction will mean that $\tilde{\gamma}_{\comps(\sigma, k)} = \tilde{\gamma}_k$ for all $k \in [K' - 1]$, and $\tilde{\gamma}_{\comps(\sigma, K')} = \tilde{\gamma}_{K'} + \sum_{i : \Lambda_i = 0}\gamma_i$. Thus  $\sum_{\ell = 1}^k\comps{\sigma, \ell} = \sum_{\ell = 1}^k \tilde{\gamma}_k > 0$ for all $k \in [K' - 1]$, and  $\sum_{\ell = 1}^{K'}\comps{\sigma, \ell} = |\gamma| > 0$ as required.
\hfill $\Box$

{\sc Proof of Proposition~\ref{prop:delay_min_DAG}}: Because the total delays are weighted averages of conditional delays, we know if the only conditional delay we are taking the average over is the minimum possible conditional delay, we will achieve the minimum total delay.  From \cref{lem:admissible_top}, we know for any admissible menu $M$, the only topological orders with positive probability are those that are admissible. 

Because the set of all permutations of CRP components is finite, the set of admissible topological orders is finite. Thus there will be some implementable topological order that achieves the minimum conditional delay (If there are some $i \in [n]$ such that $\Lambda_i = 0$, for each topological order we would also need to consider the assignment of customers classes with zero arrivals to servers that minimises delay for each topological order).

Therefore we will be able to minimise the total average delay by choosing an admissible menu $M$ that only allows for the admissible topological order that achieves the minimum conditional delay. We know that such a menu exists from \cref{lem:admissible_top}.
\hfill $\Box$
\ \\
\vspace{0.5cm}

\section{\cref{sec:waits_proof} Proofs} \label{sec:proof_waits_appen}
{\sc Proof of Lemma~\ref{lem:CRP_properties}}: There are two differences between the setup in our paper and in \cite{Afecheetal2019}: first, the constants $\gamma_i$ for the approach to heavy-traffic are allowed to be arbitrary, while in \cite{Afecheetal2019} the authors impose $\gamma_i = \Lambda_i$. Second, our setup has service classes with $\Lambda_i=0$ and hence CRP components which consist of a single service class and no servers. Despite these, the proofs for parts (i) and (ii) are identical to the proofs of parts (i) and (ii) of \cite[Lemma 3]{Afecheetal2019}.

Part (iii) of \cite[Lemma 3]{Afecheetal2019} states that $U_{\mathscr{S}_k}(M) = \mathscr{C}_k$, which in our setup should be interpreted as 
\[ U_{\mathscr{S}_k}(M) \cap \left\{ \cup_{\ell=1}^{K'} \Ccal_\ell \right\} = \mathscr{C}_k. \]
In addition, a server-less CRP component $\mathbb{C}_\kappa = ( \{i\}, \emptyset)$ consisting of a single service class $i$ is part of the set of service classes uniquely served by the set $U_{\mathscr{S}_k}(M)$ if and only if all the CRP components $k'$ such that $\mathbb{C}_\kappa$ has a directed arc to $\mathbb{C}_{k'}$ in the DAG $\Dcal=([K],\Acal)$ are included in $(\sigma(1),\ldots, \sigma(k))$. Recalling the definition of the function $\comps(\sigma, \cdot)$, this is equivalent to saying that $\comps^{-1}(\sigma, \kappa) \leq k$. 

Part (iv) follows from the definition of slack $\Delta()$ and part (iii):
\[  \Delta(\mathscr{S}_k) = \mu_{\mathscr{S}_k} - \lambda_{U_{\mathscr{S}}(M)} = \sum_{\ell = 1}^k \mu_{\Scal_\ell} - \sum_{\ell=1}^k \sum_{\kappa \in \comps(\sigma, \ell)}  \lambda_{\Ccal_\kappa} = \sum_{\ell=1}^k \sum_{\kappa \in \comps(\sigma, \ell)} \mu_{\Scal_\kappa} - \lambda_{\Ccal_\kappa} =: \epsilon \sum_{\ell = 1}^k \widetilde{\gamma}_{\comps(\sigma, \ell)} + o(\epsilon). \]

\hfill $\Box$

{\sc Proof of Lemma~\ref{lem:criticalservers}}: The first part follows from the proof of \cite[Lemma 4]{Afecheetal2019} where it is argued that if the subset $S = \{s_1, \ldots, s_\ell \}$ does not obey the condition mentioned, then
\[   \mu_S - \Lambda_{U_S(M)} > 0 , \]
and hence $\lim_{\epsilon \to 0} \frac{\epsilon}{\Delta(S)} = 0$. The second part follows from part (iv) of Lemma~\ref{lem:CRP_properties}. \hfill $\Box$
\ \\

{\sc Proof of Proposition~\ref{prop:steadystate}}: The proof of the first part of the Proposition follows exactly the same lines as \cite[Proposition 2]{Afecheetal2019} and hence we omit it. The calculations for the second part are as follows. Fix a topological ordering $\sigma \in \Tcal(\Dcal, K')$, sub-permutations $\sbb_k \in \Sigma_{\Scal_k}$, and $s = (\sbb_{\sigma(1)} || \cdots || \sbb_{\sigma(K')})$. For succinctness, define $m_k$ for $k \in \{0,1,\ldots, K'-1\}$ by
\[  m_{0} = 0 , \quad \mbox{and} \qquad m_{\ell} = m_{\ell-1} + |\Scal_{\sigma(\ell-1)}|.\]
From \eqref{eq:steadyP_k}
\begin{align*}
\pi(P(s;m)) &=  {\cal B}\, \;\prod_{\ell=1}^m \frac{1} {\Delta(s_1,\ldots, s_\ell)}    \\
&=  \Bcal \prod_{k=1}^{K'} \left( \prod_{\ell=m_{k-1}+1}^{m_{k}-1} \frac{1}{\Delta(s_1, \ldots, s_\ell)} \right) \cdot \frac{1}{\Delta(s_1,\ldots, s_{m_{k}})}.
\end{align*}
By Lemma~\ref{lem:criticalservers}, 
\[ \lim_{\epsilon \to 0} \frac{\epsilon}{\Delta(s_1,\ldots, s_{m_{k}})}  = \frac{1}{\sum_{i=1}^k \widetilde{\gamma}_{\comps(\sigma, i)}}. \]
For some $k \in [K']$, and $ m_{k-1} + 1 \leq \ell \leq m_{k}-1$, denote $S = \{ s_{m_{k-1}+1}, \ldots, s_\ell\}$. Following the outline in \cite[Lemmas 5 and 8]{Afecheetal2019}, it follows that:
\begin{align*}
    \lim_{\epsilon \to 0} \Delta(s_1,\ldots, s_{\ell}) &= \mu_{S} - \Lambda_{U_{S}(\breve{M}) } > 0.
\end{align*} 
For $\sbb_k = (s_{k}(1), \ldots. s_{k}(|\Scal_{k}|)) \in \Sigma_{\Scal_k}$, denote
\begin{equation}\label{eq:theta_k}
    \theta_k(\sbb_k) = \prod_{\ell=1}^{|\Scal_{k}|-1} \frac{1}{\mu_{\{s_k(1),\ldots, s_k(\ell)\}} - \Lambda_{U_{\{s_k(1),\ldots, s_k(\ell)\}}(\breve{M}) }}.
\end{equation}
Then,
\begin{align*}
\lim_{\epsilon \to 0}\pi(P(s;m)) &=  \lim_{\epsilon \to 0} \frac{\Bcal}{\epsilon^{K'}} \prod_{k=1}^{K'} \left( \prod_{\ell=m_{k-1}+1}^{m_{k}-1} \frac{1}{\Delta(s_1, \ldots, s_\ell)} \right) \cdot \frac{\epsilon}{\Delta(s_1,\ldots, s_{m_{k}})} \\
&= \Bcal' \left(\prod_{k=1}^{K'} \frac{1}{\sum_{i=1}^k \widetilde{\gamma}_{\comps(\sigma, i)}} \right) \left( \prod_{k=1}^{K'} \theta_k(\sbb_k) \right) \\
&= \Bcal' \cdot \mathbb{Q}(\sigma) \cdot \prod_{k=1}^{K'} \theta_k(\sbb_k),
\end{align*}
where $\Bcal' = \lim_{\epsilon \to 0} \Bcal \epsilon^{-K'}$.
\hfill $\Box$
\ \\

{\sc Proof of Lemma~\ref{lem:scaled_W_serv_perm}}: Let $s=(\sbb_{\sigma(1)} || \cdots || \sbb_{\sigma(K')}) = (s_1,\ldots, s_m) \in \Sigma_m$ be induced by topological order $\sigma \in \Tcal(\Dcal, K')$, and define $m_\ell$ for $\ell \in \{0,1,\ldots, K'-1\}$ by
\[  m_{0} = 0 , \quad \mbox{and} \qquad m_{\ell} = m_{\ell-1} + |\Scal_{\sigma(\ell-1)}|.\]
Define $j(s, i) = \min \{  \ell :  i \in U(s_1, \ldots, s_\ell)\}$, and define $\kappa$ satisfying $m_{\kappa-1}+1\leq j \leq m_{\kappa}$. Then, using Lemma~\ref{lem:WaitingTime}, we have 
\begin{align*}
 \lim_{\epsilon \to 0} \epsilon \cdot W_i(s;m) &=  \lim_{\epsilon \to 0} \sum_{\ell=j(s,i)}^m  \frac{\epsilon}{\Delta(s_1,\dots,s_\ell)} \\
 \intertext{and since each of $\lim_{\epsilon \to 0} \frac{\epsilon}{\Delta(s_1, \ldots, s_\ell)}$ exists by Lemma~\ref{lem:criticalservers},}
&=  \sum_{\ell=j(s,i)}^m \lim_{\epsilon \to 0} \frac{\epsilon}{\Delta(s_1,\dots,s_\ell)} \\
 &= \sum_{k=\kappa}^{K'} \lim_{\epsilon \to 0} \frac{\epsilon}{\Delta(s_1,\dots,s_{m_{k}})}  + \sum_{\substack{ j(s,i) \leq \ell \leq m, \\ \nexists k \ :\  \ell = m_k }} \lim_{\epsilon \to 0} \frac{\epsilon}{\Delta(s_1,\dots,s_\ell)}  \\
 &= \sum_{k = \kappa }^{K'} \frac{1}{\sum_{\ell=1}^k  \widetilde{\gamma}_{\comps(\sigma, \ell)} }.
\end{align*} 
 The last equality follows because the second term in the preceding expression is 0 by Lemma~\ref{lem:criticalservers}, and each of the terms in the first sum is precisely of the form \eqref{eq:criticalservers} in Lemma~\ref{lem:criticalservers}. The Lemma now follows by noting that $\kappa$ only depends on the CRP component $\mathbb{C}_k$ that service class $i$ belongs to and therefore so does the last expression, and $\kappa = \comps^{-1}(\sigma, k)$.
\hfill $\Box$

\section{\cref{sec:proof_matching} Proofs} \label{appn:proof_matching}

{\sc Proof of Lemma~\ref{lem:prob_transition_out}}: Let $\Scal'$ be the set of all server permutations that are not induced by any topological order. Let $s$ be a server permutation induced by some topological order $\sigma \in \Tcal(\Dcal, K')$. 

We know from flow balance that 
\begin{equation*}
    \lim_{\epsilon \to 0} \sum_{s' \in \Scal'}\sum_{b = 0}^m \pi(P(s',b)) \ge \lim_{\epsilon \to 0} \sum_{x \in P(s,m)}\pi(x)q_{ij}^0(x).
\end{equation*}
\hfill $\Box$

But \cref{prop:steadystate} tells us that 
\begin{equation*}
    \lim_{\epsilon \to 0} \sum_{s' \in \Scal'}\sum_{b = 0}^m \pi(P(s',b)) = 0.
\end{equation*}
Since $\pi(x) \in [0, 1]$ and $q_{ij}^0(x) \in [0, 1]$ for all $i \in [n]$, $j \in [m]$, and $x \in P(s,m)$, this means that 
\begin{equation*}
    \lim_{\epsilon \to 0} \sum_{x \in P(s,m)}\pi(x)q_{ij}^0(x) = 0. 
\end{equation*}

{\sc Proof of Lemma~\ref{lem:q_ij_limit}}:
Recall from \cref{def:permutationinduced} that since the permutation of servers $s$ is induced by the topological order $\sigma$, we can express $s$ as the concatenation of sub-permutations: 
\[ s = \left( \sbb_{\sigma(1)} || \sbb_{\sigma(2)}|| \cdots || \sbb_{\sigma(K')} \right) \]
with $\sbb_\kappa \in \Sigma_{\Scal_\kappa}$ denoting a permutation of the servers $\Scal_{\kappa}$ of CRP component $\mathbb{C}_\kappa$.

For $\sbb_\kappa = (s_{\kappa}(1), \ldots. s_{\kappa}(|\Scal_{\kappa}|)) \in \Sigma_{\Scal_\kappa}$, denote
\[ \theta_\kappa(\sbb_\kappa) = \prod_{\ell=1}^{|\Scal_{\kappa}|-1} \frac{1}{\mu_{\{s_{\kappa(1)},\ldots, s_\kappa(\ell)\}} - \Lambda_{U_{\{s_\kappa(1),\ldots, s_\kappa(\ell)\}}(\breve{M}) }}.\]

Also denote for $s_k \in \Sigma_k$
\begin{align*}
    H_{ij}(s_k) = & \lim_{\epsilon \to 0} \sum_{r = \hat{j}}^{|\Scal_k| - 1} \Bigg[  \left( \prod_{u = \hat{j}}^r \frac{1}{\Delta_j(s_1, \dots, s_u)} \right) \left( \prod_{\ell = r + 1}^{|\Scal_k|-1} \frac{1}{\Delta(s_1, \dots, s_\ell)} \right) \nonumber \\
    & \times \left( \frac{1}{\Delta(s_1, \dots, s_r)} - \frac{1}{\Delta_j(s_1, \dots, s_r)} \right) \Bigg] + \prod_{u = \hat{j}}^{|\Scal_k|}\frac{1}{\Delta_j (s_1, \dots, s_u)}
\end{align*}
and 
\begin{align*}
    G_{ij}(s_k) & = \lim_{\epsilon \to 0} \frac{1}{\Delta_j(s_{k}(1),\dots,s_k(|\Scal_k|)}\prod_{u = \hat{j}}^{|\Scal_k|}\frac{1}{\Delta_j (s_1, \dots, s_u)}.
\end{align*}
Finally also recall the definition of $\mathbb{Q}(\sigma)$ from \cref{eq:Q_def} as 
\begin{align*}
    \mathbb{Q}(\sigma) &= \prod_{\kappa \in [K']} \frac{1}{\sum_{\ell=1}^\kappa  \widetilde{\gamma}_{\comps(\sigma_t, \ell)}}.
\end{align*}

This lets us write $q_{ij}(P(s,m)) = \lim_{\epsilon \to 0}\htp{q}_{ij}(P(s,m))$ as 
\begin{align}
q_{ij}(P(s,m)) & = \frac{\Bcal'\lambda_i}{\pi(P(s,m))} \mathbb{Q}(\sigma) \left(\prod_{\kappa \neq k} \theta_\kappa(s_\kappa) \right) H_{ij}(s_k) \nonumber \\
& - \lim_{\epsilon \to 0} \left[ \frac{\epsilon\Bcal'\lambda_i}{\pi(P(s,m))} \left(\prod_{\kappa \neq k} \frac{1}{\sum_{\ell = 1}^\kappa \widetilde{\gamma}_{\comps(\sigma, \ell)}} \right) \left(\prod_{\kappa \neq k} \theta_\kappa(s_\kappa) \right) G_{ij}(s_k) + o(\epsilon),\right]
\end{align}
where $\Bcal' = \lim_{\epsilon \to 0} \Bcal \epsilon^{-K'}$.
\hfill $\Box$

{\sc Proof of Lemma~\ref{lem:permutation_prob}}: From \cref{prop:steadystate}, we know that 
\begin{equation}
   \lim_{\epsilon \to 0}\pi(P(s, m)) = \Bcal' \cdot \mathbb{Q}(\sigma) \prod_{k=1}^{K'} \theta_k( \sbb_k),
\end{equation}
where $\theta_k(s_k)$ is given by \cref{eq:theta_k}.

From the definition of $P_k(\sbb_k)$,we have that
\begin{align}
    \pi(P_k(\sbb_k)) = \sum_{\sigma \in \Tcal(\Dcal, K')} \ \ \sum_{\substack{  s = (\sbb_{\sigma(1)} || \sbb_{\sigma(2)} || \cdots || \sbb_k || \cdots || \sbb_{\sigma(K')})  \\ \{\sbb_k \in \Sigma_{\Scal_\kappa}\}_{\kappa \in [K']} } } \pi (P(s,m)).
\end{align}

This means that 
\begin{align}
   \lim_{\epsilon \to 0}\pi_M(P(s_k)) & = \Bcal_{M}' \sum_{\sigma \in \Tcal(\Dcal, K')} \ \ \sum_{\substack{  s = (\sbb_{\sigma(1)} || \sbb_{\sigma(2)} || \cdots || \sbb_k || \cdots || \sbb_{\sigma(K')})  \\ \{\sbb_k \in \Sigma_{\Scal_\kappa}\}_{\kappa \in [K']} } } \mathbb{Q}(\sigma) \prod_{k=1}^{K'} \theta_k( \sbb_k)  \nonumber \\
   & = \Bcal_{M}' \sum_{\sigma \in \Tcal(\Dcal, K')}  \Bigg[ \mathbb{Q}(\sigma) \sum_{\substack{  s = (\sbb_{\sigma(1)} || \sbb_{\sigma(2)} || \cdots || \sbb_k || \cdots || \sbb_{\sigma(K')})  \\ \{\sbb_k \in \Sigma_{\Scal_\kappa}\}_{\kappa \in [K']} } } \prod_{k=1}^{K'} \theta_k( \sbb_k) \Bigg]
\end{align}
Since the values of $\theta_\kappa(s_\kappa)$ are independent of each other and do not depend on $\sigma$, we can rewrite this as
\begin{align}\label{eq:prob_subperm}
   \lim_{\epsilon \to 0}\pi_M(P(s_k)) & = \Bcal_{M}' \cdot \theta_k(s_k) \left( \sum_{\sigma \in \Tcal(\Dcal, K')} \mathbb{Q}(\sigma) \right)  \prod_{\kappa \neq k} \sum_{\sbb_\kappa \in \Sigma_{\Scal_\kappa} } \theta_\kappa( \sbb_\kappa) 
\end{align}

Recall from \cref{sec:waits_proof}
\[ \left( \Bcal_M' \sum_{ \{\sbb_k \in \Sigma_{\Scal_k}\}_{k \in [K']} }   \prod_{k=1}^{K'} \theta_k(\sbb_k) \right) = \frac{1}{\sum_{\sigma \in \Tcal(\Dcal, K')} \mathbb{Q}(\sigma)}.\]
This lets us rewrite $\Bcal_M'$ as 
\[ \Bcal_M'  = \frac{1}{\left( \prod_{\kappa=1}^{K'}  \sum_{ \{\sbb_\kappa \in \Sigma_{\Scal_\kappa}\} } \theta_\kappa(\sbb_\kappa) \right)\sum_{\sigma \in \Tcal(\Dcal, K')} \mathbb{Q}(\sigma)}\]

Substituting this back into \cref{eq:prob_subperm}, we have that 
\begin{align}
   \lim_{\epsilon \to 0}\pi(P_k(\sbb_k)) & =   \frac{\theta_k(\sbb_k)}{ \sum_{\sbb_\kappa \in \Sigma_{\Scal_k} } \theta_\kappa( \sbb_\kappa) }.
\end{align}
But $\theta_k(\sbb_k)$ depend only on $\Lambda$, $\mu$, and $\breve{M}$, for all $k \in [K']$, proving the result. \hfill $\Box$

\end{appendices}

\end{document}